\documentclass[preprint,amsmath,amssymb,aip,jcp,nobibnotes]{revtex4-1}
\usepackage{amssymb,latexsym,amsmath,amsxtra}
\usepackage{amsfonts}
\usepackage{xspace}
\usepackage[final]{graphics}
\usepackage{color,psfrag} 
\usepackage[labelformat=simple]{subfig}
\usepackage{enumitem}
\usepackage{siunitx}
\usepackage[colorlinks,bookmarks,urlcolor=blue,citecolor=blue,linkcolor=blue]{hyperref}

\usepackage{bbm}

\renewcommand{\vec}[1]{\boldsymbol{#1}}
\newcommand{\vecs}[1]{\boldsymbol{#1}}
\newcommand{\paren}[1]{\left(#1\right)}
\newcommand{\brac}[1]{\left[#1\right]}
\newcommand{\E}{\mathbb E}
\newcommand{\avg}[1]{\E[#1]}
\newcommand{\from}{\leftarrow}
\newcommand{\D}[2]{\frac{d#1}{d#2}}
\newcommand{\PD}[2]{\frac{\partial#1}{\partial#2}}

\newcommand{\lap}[1]{\Delta#1}

\newcommand{\abs}[1]{\left|#1\right|}

\newcommand{\beqparen}[1]{\big(#1\big)}

\DeclareMathOperator{\prob}{Pr}
\newcommand{\comment}[1]{{#1}}

\newcommand{\DA}{D^{\textrm{A}}}
\newcommand{\DB}{D^{\textrm{B}}}
\newcommand{\DC}{D^{\textrm{C}}}

\newcommand{\bi}{\vec{i}}
\newcommand{\bj}{\vec{j}}
\newcommand{\bk}{\vec{k}}
\newcommand{\bja}{\bj^{a}}
\newcommand{\bjb}{\bj^{b}}
\newcommand{\bjc}{\bj^{c}}

\newcommand{\ve}{\vec{e}}

\newcommand{\veta}{\vecs{\eta}}

\newcommand{\vx}{\vec{x}}
\newcommand{\vy}{\vec{y}}

\newcommand{\vq}{\vec{q}}

\newcommand{\vO}{\vec{0}}

\newcommand{\rb}{r_{\textrm{b}}}

\newcommand{\vqa}{\vec{q}^{a}}
\newcommand{\vqb}{\vec{q}^{b}}
\newcommand{\vqc}{\vec{q}^{c}}

\newcommand{\ind}{\mathbbm{1}}

\def\R{\mathbb{R}}

\newcommand{\Z}{\mathbb{Z}}

\begin{document}

\title{A Convergent Reaction-Diffusion Master Equation}

\author{Samuel A. Isaacson}
\email{isaacson@math.bu.edu}
\noaffiliation
\affiliation{Department of Mathematics and Statistics\\ Boston University\\ Boston, MA 02215, USA}


\begin{abstract}
  The reaction-diffusion master equation (RDME) is a lattice
  stochastic reaction-diffusion model that has been used to study
  spatially distributed cellular processes. \comment{The RDME is often
    interpreted as an approximation to spatially-continuous models in
    which molecules move by Brownian motion and react by one of
    several mechanisms when sufficiently close. In the limit that the
    lattice spacing approaches zero, in two or more dimensions, the
    RDME has been shown to lose bimolecular reactions. The RDME is
    therefore not a convergent approximation to any
    spatially-continuous model that incorporates bimolecular
    reactions.  In this work we derive a new convergent RDME (CRDME)
    by finite volume discretization of a spatially-continuous
    stochastic reaction-diffusion model popularized by Doi.} We
  demonstrate the numerical convergence of reaction time statistics
  associated with the CRDME. 
  For sufficiently large lattice spacings or slow bimolecular reaction
  rates, we also show the reaction time statistics of the CRDME may
  be approximated by those from the RDME. 
  The original RDME may therefore be interpreted as an
  approximation to the CRDME in several asymptotic limits.
\end{abstract}

\maketitle


Computational models of biochemical systems within individual cells
have become a common tool used in studying cellular processes and
behavior~\cite{IyengarCellShape08,LauffenbergerHer22003,
  OdellSegPolNew2000,StolovitzkyP53Model2005}.
The spatially distributed nature of chemical pathways inside 
cells may, in certain cases, be more accurately modeled by including
the explicit spatial movement of molecules.  Examples of such
processes include whether signals can propagate from the plasma
membrane to nucleus~\cite{MunozGarcia:2009hd}, how cell shape can
modify information flow in signaling
networks~\cite{IyengarCellShape08}, how the variable density of
chromatin influences the search time of proteins for DNA binding
sites~\cite{IsaacsonPNAS2011}, or how different regions of cytosolic
space may separate to different chemical states~\cite{ElfIEESys04}.

One important consideration in developing these models is that at the
scale of a single cell many biochemical processes are
stochastic~\cite{ArkinPNAS1997, CollinsNatureEukNoise,
  RaserScience2004}. Stochastic reaction-diffusion models have been
used to account for the stochasticity inherent in the chemical
reaction process and the diffusion of proteins and mRNAs.  These
models approximate \emph{individual molecules} as \comment{points or
  spheres} diffusing within cells.
They are more macroscopic descriptions than quantum mechanical or
molecular dynamics models, which can resolve detailed interactions
between a few molecules on timescales of
milliseconds~\cite{ShawAntonMS2009}. They are more microscopic
descriptions than \emph{deterministic} three-dimensional
reaction-diffusion PDEs for the average concentration of each species
of molecule.

Three stochastic reaction-diffusion models that have been used to
study cellular processes are: the Doi
model~\cite{TeramotoDoiModel1967,DoiSecondQuantA,DoiSecondQuantB}, the
Smoluchowski diffusion limited reaction
model~\cite{SmoluchowskiDiffLimRx,KeizerJPhysChem82}, and the
reaction-diffusion master equation
(RDME)~\cite{GardinerRXDIFFME,GardinerHANDBOOKSTOCH,VanKampenSTOCHPROCESSINPHYS,ErbanRDMEReview,IsaacsonSJSC2006,IsaacsonRDMELims}.
In the Doi model~\cite{DoiSecondQuantA,DoiSecondQuantB} positions of
molecules are represented as points undergoing Brownian motion.
Bimolecular reactions between two molecules occur with a fixed
probability per unit time when two reactants are separated by
\emph{less} than some specified ``reaction radius''.  \comment{The
  Smoluchowski model differs by 
  representing bimolecular reactions in one of two ways; either
  occurring instantaneously (called pure absorption), or with fixed
  probability per unit time (called partial absorption), when two
  reactants' separation is \emph{exactly} the
  reaction-radius~\cite{SmoluchowskiDiffLimRx,KeizerJPhysChem82}. Two
  possible reactants are not allowed to approach closer than their
  reaction-radius. A more microscopic version of the Smoluchowski
  model approximates molecules as spheres, allowing for volume
  exclusion and collisions between non-reactive molecules.}  In each
of the models unimolecular reactions represent internal processes.
They are assumed to occur with exponentially distributed times based
on a specified reaction-rate constant.

The RDME can be interpreted as an extension of the non-spatial
chemical master equation
(CME)~\cite{GardinerRXDIFFME,McQuarrieJAppProb,GardinerHANDBOOKSTOCH,VanKampenSTOCHPROCESSINPHYS}
model for stochastic chemical kinetics.  In the RDME space is
partitioned by a mesh into a collection of voxels. The diffusion of
molecules is modeled as a continuous time random walk on the mesh,
with bimolecular reactions occurring with a fixed probability per unit
time for molecules within the same voxel.  Within each voxel molecules
are assumed well-mixed (\textit{i.e.}  uniformly distributed), and
molecules of the same species are indistinguishable.
\comment{Molecules are treated as points in the RDME, and in the
  absence of chemical reactions the continuous time random walk each
  molecule undergoes converges to the Brownian motion of a point
  particle as the mesh spacing approaches
  zero~\cite{IsaacsonSJSC2006,IsaacsonRDMENote,LodstadtRDMEFE09}.}
Mathematically, the RDME is the forward Kolmogorov equation for a
continuous-time jump Markov process.

In each of the three models the state of the chemical system is given
by stochastic processes for the number of molecules of each chemical
species and the corresponding positions of each molecule
\comment{(where in the RDME a molecule's position corresponds to the
  voxel containing that molecule)}. There are two different
mathematical formulations of each model; either coupled systems of
equations for the probability densities of having a given state at a
specified time, or equations for the evolution of the stochastic
processes themselves.  The former description leads to, possibly
infinite, coupled systems of partial integral differential equations
(PIDEs) for the Doi/Smoluchowski models, and, possibly infinite,
coupled systems of ordinary differential equations (ODEs) for the
RDME.  The high dimensionality of these equations for typical
biological systems prevents the use of standard numerical methods for
PDEs/ODEs in low-dimensions. Instead, the probability density
solutions to these equations are approximated by simulating the
underlying stochastic processes with Monte Carlo
methods~\cite{ErbanChapman2009,Bartol2008,AndrewsBrayPhysBio2004,
  DonevJCP2010,WoldeEgfrdPNAS2010,GillespieJPCHEM1977,KalosKMC75,Arjunan:2009fd}.

\comment{When the RDME is interpreted as an independent physical
  model, there is a natural, nonzero, lower bound on the lattice
  spacing in systems with bimolecular
  reactions~\cite{ElfIEESys04,IsaacsonSJSC2006,Hellander:2012jk}. In
  its most basic form, this bound arises from the physical assumption
  that the timescale for molecules to become well-mixed within a voxel
  is smaller than that for bimolecular reactions to occur. The use of
  a well-mixed bimolecular reaction mechanism within each voxel is
  then physically justified.  While the RDME can be used as an
  independent physical model, in applications it is often interpreted
  as a formal approximation of the Doi or Smoluchowski
  models~\cite{IsaacsonRDMENote,ElfPNASRates2010}.  It is important to
  stress that the RDME can probably be interpreted as a
  (non-convergent) physical approximation to a large number of
  distinct, microscopic, spatially-continuous models. For example,
  Gillespie has developed an argument for the physical validity of the
  RDME with sufficiently large lattice spacings as a (non-convergent)
  approximation to a model that involves the diffusion of hard-spheres
  which may react upon collision~\cite{GillespieComments2013}.

  In this work we \emph{interpret} the RDME as an attempt to
  approximate to the Doi model. While previous work on improving the
  accuracy of the RDME as an approximation to microscopic models has
  focused on the (point-particle) version of the Smoluchowski
  model~\cite{Hellander:2012jk,ElfPNASRates2010}, we showed
  in~\cite{IsaacsonRDMENote} that the RDME can be seen to approximate
  a version of the Doi model.  A number of studies of microscopic
  stochastic reaction-diffusion processes have been based on the Doi
  model~\cite{Naumann1993kb,Kuzovkov1988vk,Seki2007im,Klann2012ko,ErbanChapman2009,ErbanChapman2011},
  including recent work modeling the dynamics of intracellular calcium
  release~\cite{FleggCalcium2013bc}. As we discuss in the next
  section, recent studies~\cite{IsaacsonAgbanBMB2013,ErbanChapman2011}
  demonstrate that the Doi model offers comparable accuracy to the
  point-particle version of the Smoluchowski model.

  The assumption that molecules are represented by point particles is
  consistent with the use of continuous time random walks to
  approximate molecular diffusion.  In systems without bimolecular
  reactions, the RDME should converge as the lattice spacing approaches
  zero to a model of point particles that move by Brownian motion and
  may undergo linear reactions. For systems in which no reactions are
  allowed, the RDME can be shown to converge to a model for a
  collection of point-particles that undergo independent Brownian
  motions~\cite{IsaacsonRDMENote}. This model is the standard starting
  point for deriving modified spatial hopping rates in the RDME when
  incorporating new spatial transport
  methods~\cite{ElstonPeskinJTB2003,IsaacsonPNAS2011,HellanderActTransport2010}
  or other types of lattice
  discretizations~\cite{IsaacsonSJSC2006,LodstadtRDMEFE09,HellanderActTransport2010,Bayati:2011kj}.
  In the remainder, when we refer to the Smoluchowski model we will
  assume that molecules are represented by points unless stated
  otherwise.

  In practice, using the RDME to approximate either of the Doi or
  Smoluchowski models for systems including bimolecular reactions can
  be problematic. It has been shown that in the continuum limit where
  the mesh spacing in the RDME is taken to zero bimolecular reactions
  are lost~\cite{IsaacsonRDMELims, Hellander:2012jk} (in two or more
  dimensions).  That is, the time at which two molecules will react
  becomes infinite as the mesh spacing is taken to zero. This result
  is consistent with the physical lower bound on the lattice spacing,
  and as such the RDME can only provide an approximation to the Doi or
  Smoluchowski models for mesh spacings that are neither too large or
  too small~\cite{IsaacsonRDMELimsII}. The error in this approximation
  can not be made arbitrarily small; in~\cite{IsaacsonRDMELimsII} we
  found that for biologically relevant parameters values, the standard
  RDME could at best approximate within five to ten percent the
  binding time distribution for the two-molecule $\textrm{A} +
  \textrm{B} \to \varnothing$ in the Smoluchowski model (in $\R^3$). }

In this work we offer one possible method to overcome this lower limit
on the approximation error. We discretize the Doi model so as to
obtain a forward Kolmogorov equation describing a continuous-time jump
Markov process on a lattice similar to the RDME.  We seek a
discretization of this form so that we may use Monte Carlo methods to
simulate the underlying stochastic process for systems of many
reacting molecules.  The convergent RDME (CRDME) we obtain
approximates the Doi model, but retains bimolecular reactions as the
mesh spacing approaches zero.  \comment{ This is accomplished by
  decoupling the region in which two molecules can react from the
  choice of mesh. Consider the two molecule $\textrm{A} + \textrm{B}
  \to \varnothing$ reaction. In the RDME, an \textrm{A} molecule
  within a given voxel can only react with \textrm{B} molecules within
  the same voxel. The CRDME enforces that molecules with separation
  smaller than a specified reaction-radius can react, as in the Doi
  model. In the CRDME molecules are assumed well-mixed within the
  voxel containing them, and so their position is known only to the
  scale of one voxel. An \textrm{A} molecule within a given voxel can
  then react with any \textrm{B} molecule in any nearby voxel where
  the minimal distance between points within the two voxels is less
  than the reaction radius.  In this way, when the lattice spacing is
  \emph{greater} than the reaction radius, an \textrm{A} molecule
  within a specified voxel can only react with \textrm{B} molecules in
  nearest neighbor voxels (including diagonal neighbors). As the
  lattice spacing is made \emph{smaller} than the reaction-radius,
  reactions can occur between molecules that are separated by multiple
  voxels. The number of voxels separating two voxels for which a
  reaction can occur increases as the lattice spacing is decreased.
  The probability per unit time for two molecules to react in the
  CRDME is a non-increasing function of the separation between the
  voxels containing the two molecules. It decreases to zero for all
  pairs of voxels in which the minimal distance between any two points
  is larger than the reaction radius.}

For $T$ the random variable for the time at which two molecules react,
we show in two-dimensions by numerical simulation that both the
survival time distribution, $\prob \brac{T > t}$, and the mean
reaction time converge to finite values as the mesh spacing approaches
zero in the CRDME (in contrast to the RDME).  \comment{We also verify
  by comparison with Brownian dynamics
  simulations~\cite{ErbanChapman2009,ErbanChapman2011} that both the
  survival time distribution and mean reaction time converge to that
  of the Doi model.}

The new CRDME retains many of the benefits of the original RDME model,
and allows the re-use of the many extensions of the RDME that have
been developed.  Examples of these extensions include Cartesian-grid
methods for complex geometries~\cite{IsaacsonSJSC2006}; the
incorporation of drift due to potentials~\cite{IsaacsonPNAS2011,
  ElstonPeskinJTB2003}; advective velocity
fields~\cite{HellanderActTransport2010}; non-Cartesian
meshes~\cite{LodstadtRDMEFE09}; time dependent
domains~\cite{bakerYates2012vf}; GPU optimized simulation
methods~\cite{Vigelius:2012em}; adaptive mesh refinement techniques
(AMR) to improve the approximation of molecular
diffusion~\cite{Bayati:2011kj}; and multiscale couplings to more
macroscopic models~\cite{Ferm:2010cw}. 

In addition to retaining bimolecular reactions as the mesh spacing
approaches zero, we show that the CRDME is approximable by the RDME
for appropriate mesh spacings and parameter choices. \comment{It
  should be noted that the approach we take in no way invalidates the
  standard RDME as an approximation to microscopic,
  spatially-continuous models for sufficiently large lattice spacings.
  Instead, the new CRDME offers an approximation of the Doi model in
  which the approximation error can be controlled (through mesh
  refinement). The finite volume discretization approach we use also
  suggests a possible method for trying to derive CRDMEs that
  approximate the Smoluchowski model.}

\comment{We begin in the next section by reviewing several of the previous
approaches that have been used to improve the accuracy of the RDME in
approximating spatially continuous particle-based models.} In
Section~\ref{S:mathModelFormulation} we introduce the general Doi and
RDME models for the multiparticle $\textrm{A} + \textrm{B} \to
\textrm{C}$ reaction. These abstract formulations illustrate that all
bimolecular reaction terms in the Doi model correspond to identical
two-body interactions.  We next derive our new CRDME for the
two-particle annihilation reaction $\textrm{A} + \textrm{B} \to
\varnothing$ in Section~\ref{S:CRDME}.  The numerical convergence of
reaction time statistics associated with this CRDME are demonstrated
in Section~\ref{S:CRDMEConv}.  In Section~\ref{S:rdmeApproxCRDME} we
show how the RDME may be interpreted as an asymptotic approximation to
the CRDME for large mesh sizes or small reaction radii. The CRDME
allows molecules in neighboring voxels to react, suggesting the
question of where to place newly created reaction products for the
$\textrm{A} + \textrm{B} \to \textrm{C}$ reaction.  We assume
bimolecular reaction products in the Doi model are placed at the
center of mass of the two reactants.  In Section~\ref{S:crdmeRxProd}
we derive a general CRDME approximation of the Doi model for where to
place a newly created $\textrm{C}$ molecule following the bimolecular
reaction, $\textrm{A} + \textrm{B} \to \textrm{C}$. Finally, we
conclude in Section~\ref{S:crdmeSSA} by summarizing how the standard
stochastic simulation algorithm (SSA) for generating realizations of
the stochastic process described by the RDME should be modified for
the CRDME.

\comment{
\section{Previous Work}

There are several modified RDME models that have been developed to
improve the approximation of bimolecular
reactions~\cite{ErbanChapman2009,ElfPNASRates2010,Hellander:2012jk}.
To our knowledge, none of these works provide $\emph{convergent}$
approximations of a spatially-continuous stochastic reaction-diffusion
model as the lattice spacing is taken to zero.  Instead, they are each
designed to more accurately approximate one specific fixed statistic
over a range of lattice spacings that are above some critical
size.  The works of~\cite{ErbanChapman2009,Hellander:2012jk} derive
modified, lattice spacing dependent reaction rates for two molecules
within the same voxel.  In~\cite{Hellander:2012jk} the reaction rate
is modified to try to match the mean association time for the
two-molecule $\textrm{A} + \textrm{B} \to \varnothing$ reaction in the
Smoluchowski model. In~\cite{ErbanChapman2009} it is chosen to try to
recover the stationary distribution of the corresponding non-spatial,
well-mixed system.  These procedures only work for sufficiently large
lattice spacings~\cite{Hellander:2012jk}, for example a factor of
$\pi$ times the reaction-radius in three-dimensions for the method
of~\cite{Hellander:2012jk}.

The approach of~\cite{ElfPNASRates2010} involves several modifications
to the RDME.  The authors begin by considering the reaction
$\textrm{A} + \textrm{B} \leftrightarrows \textrm{C}$ for a system in
which there is one stationary $\textrm{A}$ molecule located at the
origin, and one $\textrm{B}$ molecule that may diffuse within a
concentric sphere about the $\textrm{A}$ molecule.  (The $\textrm{A}$
molecule is assumed to be a sphere.) The corresponding
radially-symmetric Smoluchowski model for the separation of the
$\textrm{B}$ molecule from the $\textrm{A}$ molecule is then
discretized into the form of a master equation, where the $\textrm{B}$
molecule hops on a radial mesh bounded by the surface of the
$\textrm{A}$ molecule and the outer domain boundary.  From this
discretization, the authors determine an analytic, lattice spacing
dependent bimolecular reaction rate when the $\textrm{B}$ molecule is
in the mesh voxel bordering the $\textrm{A}$ molecule. For a given
mesh spacing, this rate is chosen so that the mean equilibration time
in the lattice model is the same as in the Smoluchowski
model~\cite{ElfPNASRates2010}. The authors then adapt their reaction
rates to allow both molecules to move on a Cartesian grid lattice.
To improve the accuracy of the RDME in approximating the Smoluchowski
model, in~\cite{ElfPNASRates2010} reactions are allowed between
molecules within nearest neighbor voxels along each coordinate
direction.  For example, in two dimensions this leads to a five-point
reaction stencil.  Within this stencil, bimolecular reactions are
modeled as well-mixed, and occur with a fixed rate based on volume
corrections to the rate derived for the spherically symmetric model.

While the method of~\cite{ElfPNASRates2010} makes use of reactions
between molecules in different voxels it is important to note it is
fundamentally different than the CRDME we derive in
Section~\ref{S:CRDME}.  First, the method of~\cite{ElfPNASRates2010}
is designed to accurately approximate one particular statistic of the
Smoluchowski model, the mean equilibration time, for a \emph{large
  range} of lattice spacings that are bigger than some multiple of the
reaction radius. (The simulations in~\cite{ElfPNASRates2010} are
restricted to lattice spacings greater than twice the
reaction-radius.)  In contrast, our method is designed so that the
solution to the CRDME \emph{converges} to the solution of the Doi
model as the lattice spacing is taken to zero (and hence should be
able to approximate \emph{any} statistic of that model to arbitrary
accuracy, for sufficiently small lattice spacings). The method
of~\cite{ElfPNASRates2010} allows a molecule to react with other
molecules in a \emph{fixed} number of neighboring voxels as the
lattice spacing is changed. The volume of the region about one
molecule in which a reaction can occur will therefore approach zero as
the lattice spacing approaches zero. As in the standard RDME, this may
cause a loss of bimolecular reactions should one try to take the
lattice spacing to zero (which is \emph{not} the goal
of~\cite{ElfPNASRates2010}).  The CRDME decouples the region in which
two molecules may react from the lattice, as described in the
introduction. For lattice spacings \emph{greater} than the
reaction-radius an \textrm{A} molecule can react with \textrm{B}
molecules in any nearest-neighbor voxel (including diagonal
neighbors).  For lattice spacings \emph{smaller} than the reaction
radius, an \textrm{A} molecule can react with \textrm{B} molecules in
any voxels for which the minimum separation between any point in the
\textrm{A} molecule voxel and any point in the \textrm{B} molecule
voxel is smaller than the reaction radius. As such, the number of
voxels separating two voxels for which a reaction can occur increases
as the lattice spacing is decreased.  The effective volume in which a
reaction between two molecules is allowed approaches that of the Doi
model as the lattice spacing is decreased to zero.  Moreover, the
probability per unit time two molecules will react is a function of
the separation of the voxels that contain them, and is not constant
across all voxel pairs for which a reaction can occur (in contrast
to~\cite{ElfPNASRates2010}).

It should also be noted that the version of the Smoluchowski model
used in~\cite{ElfPNASRates2010,Hellander:2012jk} only takes into
account the physical size of molecules in the context of bimolecular
reactions. Volume exclusion due to molecule size is not accounted for
in the diffusive hopping rates (which are chosen to recover the
Brownian motion of point particles).  \emph{With this approximation},
the Doi model should provide a similar level of physical accuracy to
the Smoluchowski model.  For example, it was shown
in~\cite{IsaacsonAgbanBMB2013} that if molecules react instantly upon
reaching a fixed separation in the Smoluchowski model, as in the
popular Brownian dynamics simulator
Smoldyn~\cite{AndrewsBrayPhysBio2004}, then the solution to the Doi
model converges to the solution of the Smoluchowski model as the
probability per unit time for two molecules to react when separated by
less than a reaction-radius is increased to infinity. It was also
shown in~\cite{ErbanChapman2011} that the standard statistics one
might wish to match in the Smoluchowski model, such as specified
geminate recombination probabilities, can be captured by an
appropriate choice of reaction parameters in the Doi model.}

\section{General Doi and RDME models} \label{S:mathModelFormulation}
We first illustrate how the bimolecular reaction $\textrm{A} +
\textrm{B} \to \textrm{C}$ would be described by the continuum Doi
model and the standard lattice RDME model in $\R^d$.  In the Doi
model, bimolecular reactions are characterized by two parameters; the
separation at which molecules may begin to react, $\rb$, and the
probability per unit time the molecules react when within this
separation, $\lambda$.  When a molecule of species $\textrm{A}$ and a
molecule of species $\textrm{B}$ react we assume the $\textrm{C}$
molecule they produce is placed midway between them.  Note the
important point that in each of the models molecules are modeled as
points.

We now formulate the Doi model as an infinite coupled system of
partial integral differential equations (PIDEs).  Let $\vqa_{l} \in
\R^d$ denote the position of the $l$th molecule of species
$\textrm{A}$ when the total number of molecules of species
$\textrm{A}$ is $a$. The state vector of the species $\textrm{A}$
molecules is then given by $\vqa = (\vqa_{1},\dots,\vqa_{a}) \in
\R^{da}$.  Define
$\vqb$ and $\vqc$ similarly.
We denote by $f^{(a,b,c)}(\vqa,\vqb,\vqc,t)$ the probability density
for there to be $a$ molecules of species $\textrm{A}$, $b$ molecules
of species $\textrm{B}$, and $c$ molecules of species $\textrm{C}$ at
time $t$ located at the positions $\vqa$, $\vqb$, and $\vqc$.
Molecules of the same species are assumed indistinguishable.  The
evolution of $f^{(a,b,c)}$ is given by
\begin{equation} \label{eq:doiFormEvolEq}
    \PD{f^{(a,b,c)}}{t}\paren{\vqa,\vqb,\vqc,t} = \paren{L + R} f^{(a,b,c)}\paren{\vqa,\vqb,\vqc,t}.
\end{equation}
Note, with the subsequent definitions of the operators $L$ and $R$
this will give a coupled system of PIDEs over all possible values of
$(a,b,c)$. 
More general systems that allow unbounded production of certain
species would result in an infinite number of coupled PIDEs.  The
diffusion operator, $L$, is defined by
\begin{equation} \label{eq:doiFormLapOp}
  L f^{(a,b,c)} = 
  \brac{D^{\textrm{A}} \sum_{l=1}^{a} \Delta_{l}^{a}
  + D^{\textrm{B}} \sum_{m=1}^{b} \Delta_{m}^{b} + D^{\textrm{C}} \sum_{n=1}^{c} \Delta_{n}^{c}}  f^{(a,b,c)},
\end{equation}
where $\Delta_{l}^{a}$ denotes the Laplacian in the coordinate
$\vqa_{l}$ and $D^{\textrm{A}}$ the diffusion constant of species
$\textrm{A}$. $D^{\textrm{B}}$, $D^{\textrm{C}}$, $\lap_{m}^{b}$, and
$\lap_{n}^{c}$ are defined similarly.

To define the reaction operator, $R$, we introduce notations for
removing or adding a specific molecule to the state $\vqa$.  Let
\begin{align*}
  \vqa \setminus \vqa_{l} &= \paren{\vqa_{1},\dots,\vqa_{l-1},\vqa_{l+1},\dots,\vqa_{a}}, \\   \vqa \cup \vq &= \paren{\vqa_{1},\dots,\vqa_{a},\vq}.
\end{align*}
$\vqa \setminus \vq$ will denote $\vqa$ with any one component with
the value $\vq$ removed. 
Denote by $\ind_{\brac{0,\rb}}(r)$ the indicator function of the
interval $\brac{0,\rb}$. 
The Doi reaction operator, $R$, is then
\begin{multline} \label{eq:doiRxOp}
    (R f^{(a,b,c)})\paren{\vqa,\vqb,\vqc,t} = 
    \lambda \Bigg[ \sum_{l=1}^{c} \\
    \int_{R^{d}} \int_{R^{d}} \delta \paren{\frac{\vq + \vq'}{2} - \vqc_l}
    \ind_{\brac{0,\rb}} \paren{\abs{\vq - \vq'}} 
    f^{(a+1,b+1,c-1)} \paren{\vqa \cup \vq, \vqb \cup \vq', \vqc \setminus \vqc_{l}, t} \, d\vq \, d \vq' \\
    - \sum_{l=1}^a \sum_{l'=1}^b 
    \ind_{\brac{0,\rb}} \big( \big| \vqa_l-\vqb_{l'} \big| \big) f^{(a,b,c)}(\vqa,\vqb,\vqc,t) \Bigg].
\end{multline}
Let $B_{l}^{c} = \{ \vq \in \R^{d} |
 \abs{\vq-\vqc_{l}} \leq \rb/2\}$ label the set of points a reactant
could be at to produce a molecule of species $\textrm{C}$ at $\vqc_l$.
Then~\eqref{eq:doiRxOp} simplifies to
\begin{multline*} 
  (R f^{(a,b,c)})\paren{\vqa,\vqb,\vqc,t} = 
  \lambda \Bigg[ 2^d \sum_{l=1}^{c} 
  \int_{B_{l}^{C}} f^{(a+1,b+1,c-1)} \paren{\vqa \cup \vq, \vqb \cup \paren{2 \vqc_{l} - \vq}, \vqc \setminus \vqc_{l}, t} \, d\vq \\
  - \sum_{l=1}^a \sum_{l'=1}^b 
  \ind_{\brac{0,\rb}} \big( \big| \vqa_l-\vqb_{l'} \big| \big) f^{(a,b,c)}(\vqa,\vqb,\vqc,t) \Bigg].
\end{multline*}

We now describe the RDME, in a form we derived
in~\cite{IsaacsonRDMENote} that has the advantage of representing a
chemical system's state in a similar manner to the Doi model. Using
this representation allows for easier comparison of the RDME and Doi
models.  Partition $\R^{d}$ into a Cartesian lattice of voxels with
width $h$ and hypervolume $h^{d}$.  When in the \emph{same} voxel, an
$\textrm{A}$ and $\textrm{B}$ molecule may react with probability per
unit time $k / h^d$.  Here $k$ represents the macroscopic bimolecular
reaction-rate constant of the reaction $\textrm{A} + \textrm{B} \to
\textrm{C}$, with units of hypervolume per unit time.  Let $\bja_l \in
\Z^d$ denote the multi-index of the voxel centered at $h \bja_l$ that
contains the $l$th molecule of species $\textrm{A}$ when there are $a$
molecules of species $\textrm{A}$.  The position of the molecule is
assumed to be uniformly distributed (\textit{i.e.}  well-mixed) within
this voxel. Let $\bja = (\bja_1,\dots,\bja_a)$ denote the state vector
for the voxels containing the $a$ molecules of species $\textrm{A}$.
Define $\bjb$ and $\bjc$ similarly, and let
$F_{h}^{(a,b,c)}(\bja,\bjb,\bjc,t)$ denote the probability that there
are $(a,b,c)$ molecules of species $\textrm{A}$, $\textrm{B}$, and
$\textrm{C}$ at time $t$ in the voxels given by $\bja$, $\bjb$, and
$\bjc$. The RDME is then the coupled system of ODEs over all possible
values for $a$, $b$, $c$, $\bja$, $\bjb$, and $\bjc$,
\begin{equation} \label{eq:RDMEPartEq}
  \D{F_{h}^{(a,b,c)}}{t}\paren{\bja,\bjb,\bjc,t} = \paren{L_{h} + R_{h}} F_{h}^{(a,b,c)}\paren{\bja,\bjb,\bjc,t},
\end{equation}
where $L_h$ is a discretized approximation to $L$ given by
\begin{equation*}
  L_{h} F_{h}^{(a,b,c)}\paren{\bja,\bjb,\bjc,t} = 
  \paren{D^{\textrm{A}} \Delta_{h}^{a} + D^{\textrm{B}} \Delta_{h}^{b} + D^{\textrm{C}} \Delta_{h}^{c}} F_{h}^{(a,b,c)}\paren{\bja,\bjb,\bjc,t}.
\end{equation*}

Here $\Delta_{h}^{a}$ denotes the standard $da$-dimensional discrete
Laplacian acting in the $\bja$ coordinate. We define the ``standard''
$d$-dimensional discrete Laplacian acting on a mesh function,
$f(\bj)$ on $\Z^d$, by 
\begin{equation} \label{eq:lapDef}
  \Delta_{h} f(\bj) = \frac{1}{h^2} \sum_{k=1}^d \brac{ f(\bj + \ve_k)
    + f(\bj - \ve_k) - 2 f(\bj) },
\end{equation}
where $\ve_k$ denotes a unit vector along the $k$th coordinate axis of
$\R^d$.

The reaction operator, $R_h$, is given by 
\begin{multline} \label{eq:RDMEPartEqRxOp}
  \paren{ R_{h} F_{h}^{\paren{a,b,c}} } \beqparen{\bja,\bjb,\bjc,t} = 
  \frac{k}{h^d} \bigg[  \sum_{l=1}^{c}  F_{h}^{\paren{a+1,b+1,c-1}}
  \beqparen{\bja \cup \bjc_{l}, \bjb \cup \bjc_{l}, \bjc \setminus \bjc_{l},t} \\
  - \sum_{l=1}^{a} \sum_{m=1}^{b} \delta_{h} \beqparen{\bja_{l} - \bjb_{m}} \,
  F_{h}^{\paren{a,b,c}} \beqparen{\bja,\bjb,\bjc,t} \bigg],
\end{multline}
where $\delta_h( \bja_l - \bjb_m )$ denotes the Kronecker delta
function equal to one when $\bja_l = \bjb_m$ and zero otherwise.  Note
that these equations are a formal discrete approximation to the Doi
model, where two molecules may now react when in the same voxel with
rate $k / h^d$.

We have shown that the RDME~\eqref{eq:RDMEPartEq} may be interpreted
as a \emph{formal} approximation to both Doi-like and Smoluchowski
models~\cite{IsaacsonRDMENote, IsaacsonRDMELims,IsaacsonRDMELimsII}.
We proved in~\cite{IsaacsonRDMELims}, and showed numerically
in~\cite{IsaacsonRDMELimsII}, that the RDME loses bimolecular
reactions as $h \to 0$.  This was demonstrated rigorously for $d=3$,
where the time for two molecules two react was shown to diverge like
$h^{-1}$. A simple modification of the argument
in~\cite{IsaacsonRDMELims} shows bimolecular reactions are lost for
all $d > 1$, with a divergence like $\ln(h)$ for $d=2$, and like
$h^{-d+2}$ for $d>2$.  More recently asymptotic expansions were used
in~\cite{Hellander:2012jk} to show the mean reaction time becomes
infinite with the preceding rates as $h \to 0$ (for $d=2$ or $d=3$).
This loss of reaction occurs because molecules are modeled by points,
and as $h \to 0$ each voxel of the mesh shrinks to a point.  Since two
molecules must be in the same voxel to react, and in two or more
dimensions two points can not find each other by diffusion,
bimolecular reactions will never occur.

It should be noted that while the RDME loses bimolecular reactions in
the limit that $h \to 0$, we have shown that the solution to the RDME,
for \emph{fixed} values of $h$, gives an asymptotic approximation for
small $\rb$ to the solution of the Smoluchowski
model~\cite{IsaacsonRDMELims,IsaacsonRDMELimsII}. 
How accurate this approximation can be made is dependent on domain
geometry and the parameters of the underlying chemical
system~\cite{IsaacsonRDMELimsII}. In particular, $h$ must be chosen
sufficiently large that reactions within a voxel can be approximated
by a well-mixed reaction with rate $k / h^d$, while chosen
sufficiently small that the diffusion of the molecules is
well-approximated by a continuous time random walk on the
mesh~\cite{ElfIEESys04,IsaacsonRDMELims}.

In the next section we develop a new convergent RDME (CRDME) to
overcome these limitations of the standard RDME~\eqref{eq:RDMEPartEq}.

\section{A convergent RDME (CRDME)} \label{S:CRDME}

To construct a convergent RDME (CRDME) we use a finite volume
discretization of the Doi PIDEs~\eqref{eq:doiFormEvolEq}. For brevity
we illustrate our approach on a simplified version
of~\eqref{eq:doiFormEvolEq} when there is only one molecule of
$\textrm{A}$ and one molecule of $\textrm{B}$ in the system which may
undergo the annihilation reaction $\textrm{A} + \textrm{B} \to
\varnothing$. The approach we describe can be extended to general
multi-particle systems as bimolecular reactions in the Doi model only
involve multiple two-particle interactions of the same form,
see~\eqref{eq:doiFormEvolEq}.

For now we work in $d$-dimensional free-space, $\R^d$ (for most
biological models $d=2$ or $d=3$). Denote by $\vx \in \R^d$ the
position of the molecule of species $\textrm{A}$ and by $\vy \in \R^d$
the position of the molecule of species $\textrm{B}$.  In the Doi
model these molecules diffuse independently, and may react with
probability per unit time $\lambda$ when within a separation $\rb$.
($\rb$ is usually called the reaction-radius.)  We let $\mathcal{R} =
\{ (\vx,\vy) \mid \abs{\vx-\vy} < \rb \}$, and denote the indicator
function of this set by $\ind_{\mathcal{R}}(\abs{\vx - \vy})$.  The
diffusion constants of the two molecules will be given by $\DA$ and
$\DB$ respectively.

Finally, we denote by $p(\vx,\vy,t)$ the probability density the two
molecules have not reacted and are at the positions $\vx$ and $\vy$ at
time $t$.  Then~\eqref{eq:doiFormEvolEq} reduces to
\begin{equation} \label{eq:DoiSimpEq}
  \PD{p}{t}(\vx,\vy,t) = (D^{\textrm{A}} \lap_{\vx} + D^{\textrm{B}} \lap_{\vy}) p(\vx,\vy,t) 
  - \lambda \ind_{\mathcal{R}}(\abs{\vx-\vy}) p(\vx,\vy,t).
\end{equation}
(Here we have dropped the equation for the state $a=0$, $b=0$.)  

We now show how to construct a new type of RDME by discretization of
this equation. Note, while~\eqref{eq:DoiSimpEq} can be solved
analytically by switching to the separation coordinate, $\vx - \vy$,
such approaches will not work for more general chemical systems, such
as~\eqref{eq:doiFormEvolEq}.  For this reason, we illustrate our CRDME
ideas on~\eqref{eq:DoiSimpEq}. We later show in
Section~\ref{S:crdmeRxProd} how the RDME reaction
operator~\eqref{eq:RDMEPartEqRxOp} for the general
$\textrm{A}+\textrm{B} \to \textrm{C}$ reaction is modified in the
CRDME (see~\eqref{eq:genCRDMERxOp}).

\begin{figure}[t]
  \centering
  \scalebox{.4}{\includegraphics{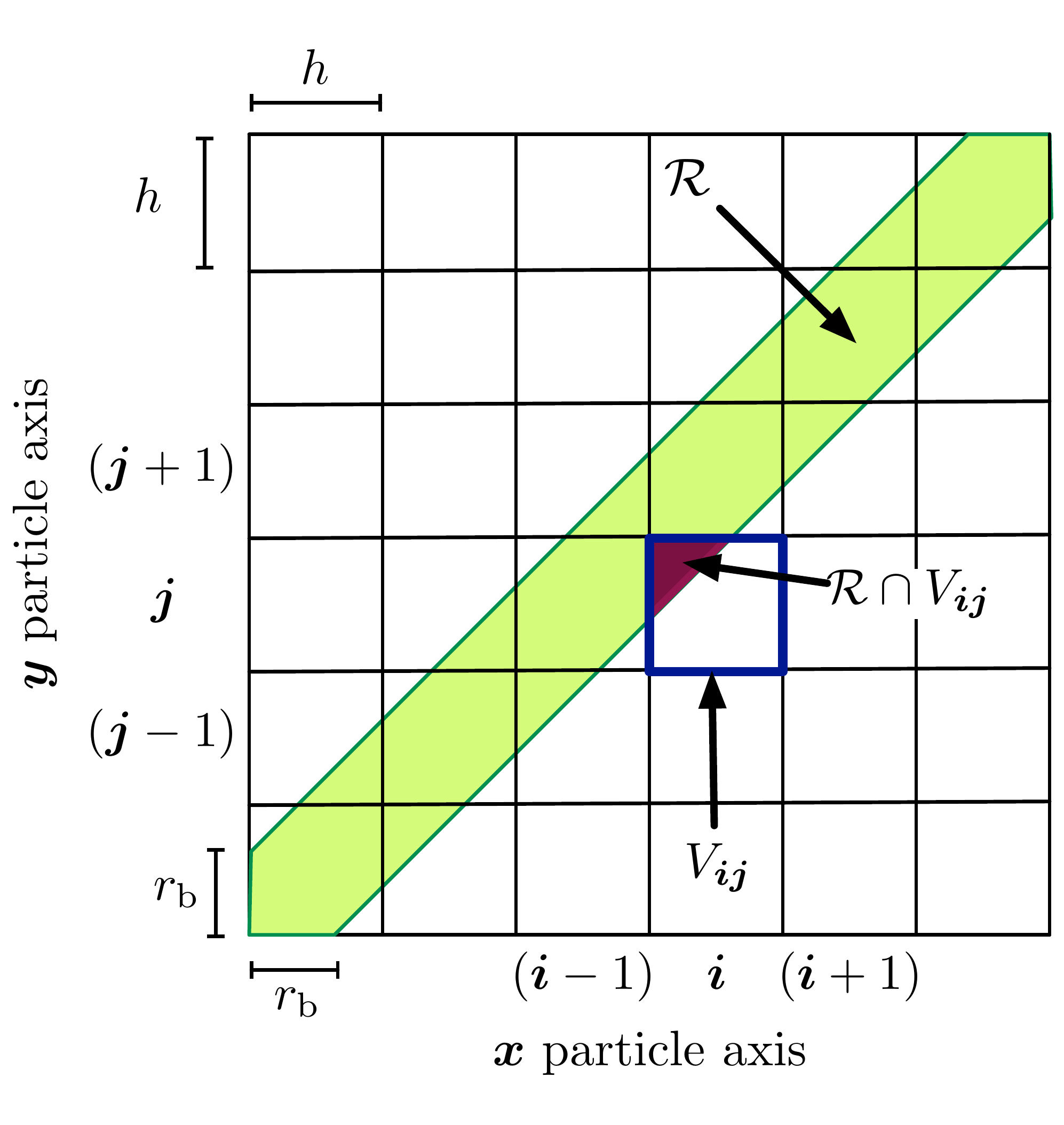}}  
  \caption{ \small \comment{Effective two-dimensional lattice when each
      particle is in $\R$. The green region corresponds to the set
      $\mathcal{R}$ of $(\vx,\vy)$ values where the two molecules can
      react. $V_{\bi \bj}$ labels the interior of the square with the
      blue boundary.  It corresponds to the set of possible
      $(\vx,\vy)$ values when the $\textrm{A}$ molecule with position
      $\vx$ is randomly distributed within $\brac{h(\bi -
        \frac{1}{2}),h(\bi + \frac{1}{2})}$, and the $\textrm{B}$
      molecule with position $\vy$ is randomly distributed within
      $\brac{h(\bj - \frac{1}{2}),h(\bj + \frac{1}{2})}$.  The maroon
      region, $\mathcal{R} \cap V_{\bi \bj}$, labels the subset of
      possible particle pair positions in $V_{\bi \bj}$ where a
      reaction can occur.}}
  \label{fig:geomDefs}
\end{figure}
For $\bi \in \Z^{d}$ and $\bj \in \Z^d$ we partition $\R^{2d}$ into a
Cartesian grid of hypercubes, labeled by $V_{\bi \bj}$. We assume
$V_{\bi \bj}$ has coordinate-axis aligned edges with length $h$ and
center $(\vx_{\bi},\vy_{\bj}) = (\bi h, \bj h)$.  Denote the
hypervolume of a set, $S$, by $\abs{S}$.  For example, the hypervolume
of the hypercube $V_{\bi \bj}$ is $\abs{V_{\bi \bj}} = h^{2d}$.
\comment{In Fig.~\ref{fig:geomDefs} we illustrate the various geometric
  quantities we will need in discretizing the Doi
  model~\eqref{eq:DoiSimpEq} when the molecules are in one-dimension
  ($d=1$).}

\comment{We make the standard finite volume method
  \emph{approximation}~\cite{DouglasFVConv:2003go} that the
  probability density, $p(\vx,\vy,t)$, is constant within each
  hypercube, $V_{\bi \bj}$. The probability the two molecules are in
  $V_{\bi \bj}$ at time $t$ is then given by $P_{\bi,\bj}(t) = p(\vx,
  \vy, t) \abs{V_{\bi \bj}}$ (for any point $(\vx,\vy) \in V_{\bi
    \bj}$).  Let $V_{\bi}$ denote the $d$-dimensional hypercube with
  sides of length $h$ centered at $\bi h$. With this definition we may
  write $V_{\bi \bj} = V_{\bi} \times V_{\bj}$.  $P_{\bi,\bj}(t)$ then
  gives the probability that the particle of species $\textrm{A}$ is
  in $V_{\bi}$, and the particle of species $\textrm{B}$ is in
  $V_{\bj}$, at time $t$. The approximation that $p(\vx,\vy,t)$ is
  constant within $V_{\bi \bj}$ is equivalent to assuming the two
  molecules are \emph{well-mixed} within $V_{\bi}$ and $V_{\bj}$
  respectively.}

Using this assumption we construct a finite volume discretization
of~\eqref{eq:DoiSimpEq} by integrating both sides
of~\eqref{eq:DoiSimpEq} over the hypercube $V_{\bi \bj}$. As we did
in~\cite{IsaacsonSJSC2006}, we make the standard finite volume
approximations for the integrals involving $\lap_{\vx} p$ and
$\lap_{\vy} p$ to obtain discrete Laplacians given
by~\eqref{eq:lapDef} in the $\vx$ and $\vy$ coordinates. The reaction
term is approximated by
\begin{align*} 
  \lambda \int_{V_{\bi \bj}} \ind_{\mathcal{R}}(\abs{\vx-\vy}) 
  p(\vx,\vy,t) \, d
  \vx \, d \vy 
   &\approx\frac{\lambda}{\abs{V_{\bi \bj}}} P_{\bi,
      \bj}(t) \int_{V_{\bi \bj}} \ind_{\mathcal{R}}(\abs{\vx-\vy})
    \, d \vx \, d \vy\\
    &= \frac{\lambda \abs{\mathcal{R} \cap V_{\bi \bj}}}{\abs{V_{\bi \bj}}} P_{\bi,
      \bj}(t).
\end{align*}
\comment{Let $\phi_{\bi \bj} = \abs{\mathcal{R} \cap V_{\bi \bj}}
  {\abs{V_{\bi \bj}}}^{-1}$ label the fraction of the total volume in
  $V_{\bi \bj}$ where a bimolecular reaction is possible (the maroon
  region in Fig.~\ref{fig:geomDefs}).  $\phi_{\bi \bj}$ is the
  probability that when the $\textrm{A}$ molecule is well-mixed in
  $V_{\bi}$ and the \textrm{B} molecule is well-mixed in $V_{\bj}$
  they are close enough to be able to react. } Our discretization then
represents a new RDME for the two-particle system, given by the
coupled system of ODEs over all values of $(\bi,\bj) \in \Z^{2d}$,
\begin{equation} \label{eq:RDMENewEq} \D{P_{\bi,\bj}}{t}(t) = L_h P_{\bi,
    \bj}(t) - \lambda \phi_{\bi \bj}  P_{\bi,\bj}(t).
\end{equation}
Here $L_h = (\DA L_{h}^{\textrm{A}} + \DB L_{h}^{\textrm{B}})$, with
$L_{h}^{\textrm{A}}$ and $L_{h}^{\textrm{B}}$ denoting standard
$d$-dimensional discrete Laplacians~\eqref{eq:lapDef} (in the $\bi$
and $\bj$ coordinates respectively). By choosing an appropriate
discretization we have obtained an equation that has the form of the
forward Kolomogorov equation for a continuous-time jump Markov
process. We may therefore interpret the coefficients
in~\eqref{eq:RDMENewEq} as transition rates, also called propensities,
between states $(\bi,\bj)$ of the stochastic process.  We subsequently
refer to this equation as the CRDME.  In the CRDME diffusion is
handled in exactly the same manner as for the RDME.
In contrast, the reaction mechanism in~\eqref{eq:RDMENewEq} now allows
molecules in distinct voxels, $V_{\bi}$ and $V_{\bj}$, to react with a
potentially non-zero probability per unit time, $\lambda \phi_{\bi
  \bj}$.

As discussed earlier, the RDME is only physically valid when $h$ is
chosen sufficiently large that the timescale for two molecules to
become uniformly distributed within a voxel by diffusion is much
faster than that for a well-mixed bimolecular reaction to occur
between them. By allowing molecules to react when in nearby voxels,
our CRDME provides a correction when $h$ is sufficiently small that
this condition is violated. For $\rb > h$ molecules may potentially
react when separated by multiple voxels, with the number of voxels
apart two molecules can be and still react increasing as $h \to 0$.

\section{Numerical convergence of the CRDME} \label{S:CRDMEConv}
We now demonstrate that in two-dimensions ($d = 2$) the survival time
distribution and the mean reaction time for the
CRDME~\eqref{eq:RDMENewEq} converge to finite values as $h \to 0$. We
show that in the corresponding RDME model the survival time and mean
reaction time diverge to $\infty$ as $h \to 0$. In contrast, in the
opposite limit that $\rb / h \to 0$ we demonstrate that the mean
reaction time in the RDME approaches that of the CRDME.
As~\eqref{eq:DoiSimpEq} is a PDE in four-dimensions, we do not
directly solve the corresponding system of ODEs given by the
CRDME~\eqref{eq:RDMENewEq}.  Instead, we simulate the corresponding
stochastic jump process for the molecule's motion and reaction by the
well-known exact stochastic simulation algorithm (SSA) (also known as
the Gillespie method~\cite{GillespieJPCHEM1977} or kinetic Monte
Carlo~\cite{KalosKMC75}).

We assume each molecule moves within a square with sides of length
$L$, $\Omega = \brac{0,L} \times \brac{0,L}$, with zero Neumann
boundary conditions.  These boundary conditions are enforced by
setting the transition rate for a molecule to hop from a given mesh
voxel outside the domain to zero. For a specified number of mesh
voxels, $N$, we discretize $\Omega$ into a Cartesian grid of squares
with sides of length $h = L / N$. We assume $\DA = \DB = D$. Unless
otherwise specified, all spatial units will be micrometers, with time
in units of seconds. The SSA-based simulation algorithm can be
summarized as:
\begin{enumerate}[itemsep=0pt]
\item Specify $D$, $L$, $N$, $\lambda$, and $\rb$ as input.
\item Calculate $\phi_{\vO \bj}$. (See
  Appendix~\ref{S:reactTransitionRates}.)
\item We assume the molecules are well-mixed at time $t=0$. That is,
  the initial position of each molecule is sampled from a uniform
  distribution among all voxels of the mesh.
\item Sample a time and direction of the next spatial hop by one
  of the molecules.  (From~\eqref{eq:RDMENewEq} each molecule may hop
  to a neighbor in the $x$ or $y$ direction with probability per unit
  time $D / h^2$.)
\item Assuming the molecules are at $\bi$ and $\bj$, if $\phi_{\bi
    \bj} \neq 0$ sample the next reaction time using the transition
  rate $\lambda \phi_{\bi \bj}$.
\item Select the smaller of the hopping and reaction times, and
  execute that event. Update the current time to the time of the
  event.
\item If a reaction occurs the simulation ends. If a spatial hop
  occurs, return to 4.
\end{enumerate}

For all simulations we chose $L = .2 \, \mu \textrm{m}$. With this
choice the domain could be interpreted as small patch of membrane
within a cell. A diffusion constant of $D = 10 \, \mu \textrm{m}^2
\textrm{s}^{-1}$ was used for each molecule. The reaction radius,
$\rb$, was chosen to be $1 \, \textrm{nm}$.  While physical reaction
radii are generally not measured experimentally, this choice falls
between the measured width of the LexA DNA binding potential ($\approx
5$ angstroms~\cite{KuhnerLexADNABond}) and the $5 \, \textrm{nm}$
reaction radius used for interacting membrane proteins
in~\cite{Dushek2011}.

We also simulated the stochastic process described by the
corresponding RDME model.  The bimolecular reaction rate was chosen to
be $k = \lambda \pi \rb^2$ to illustrate how the RDME approximates the
CRDME as $\rb / h \to 0$, but diverges as $\rb /h \to \infty$. Our
motivation for this choice is explained in the next section.  The
simulation algorithm was identical to that just described, except that
step 2 was removed and step 5 modified so that two molecules could
only react when within the same voxel (with probability per unit time
$k = \lambda \pi \rb^2 / h^2$).

\begin{figure*}
  \centering
  \subfloat[]{
    \scalebox{.54}{\includegraphics{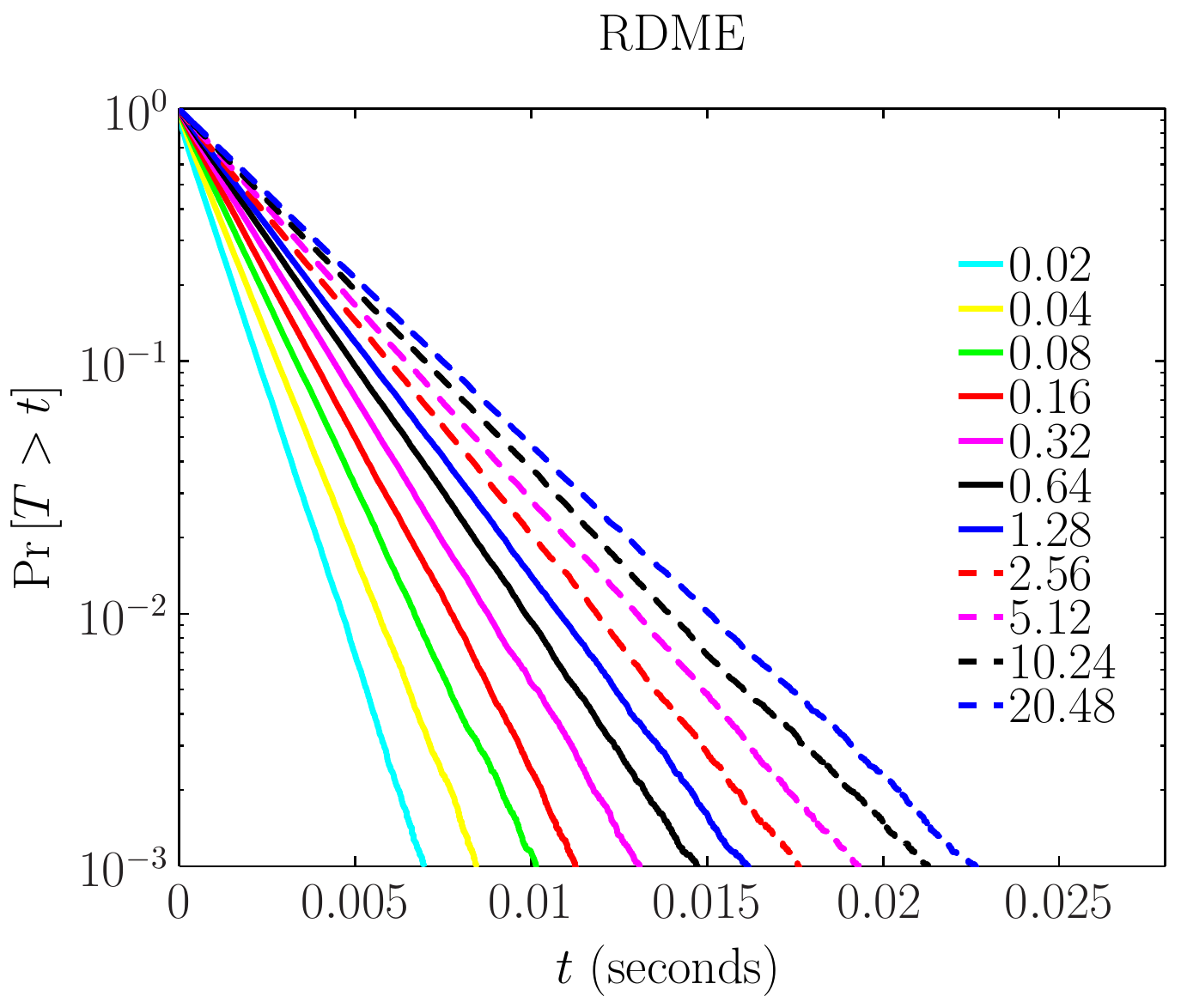}}  
  }
  \hspace{20pt}
  \subfloat[]{
    \scalebox{.54}{\includegraphics{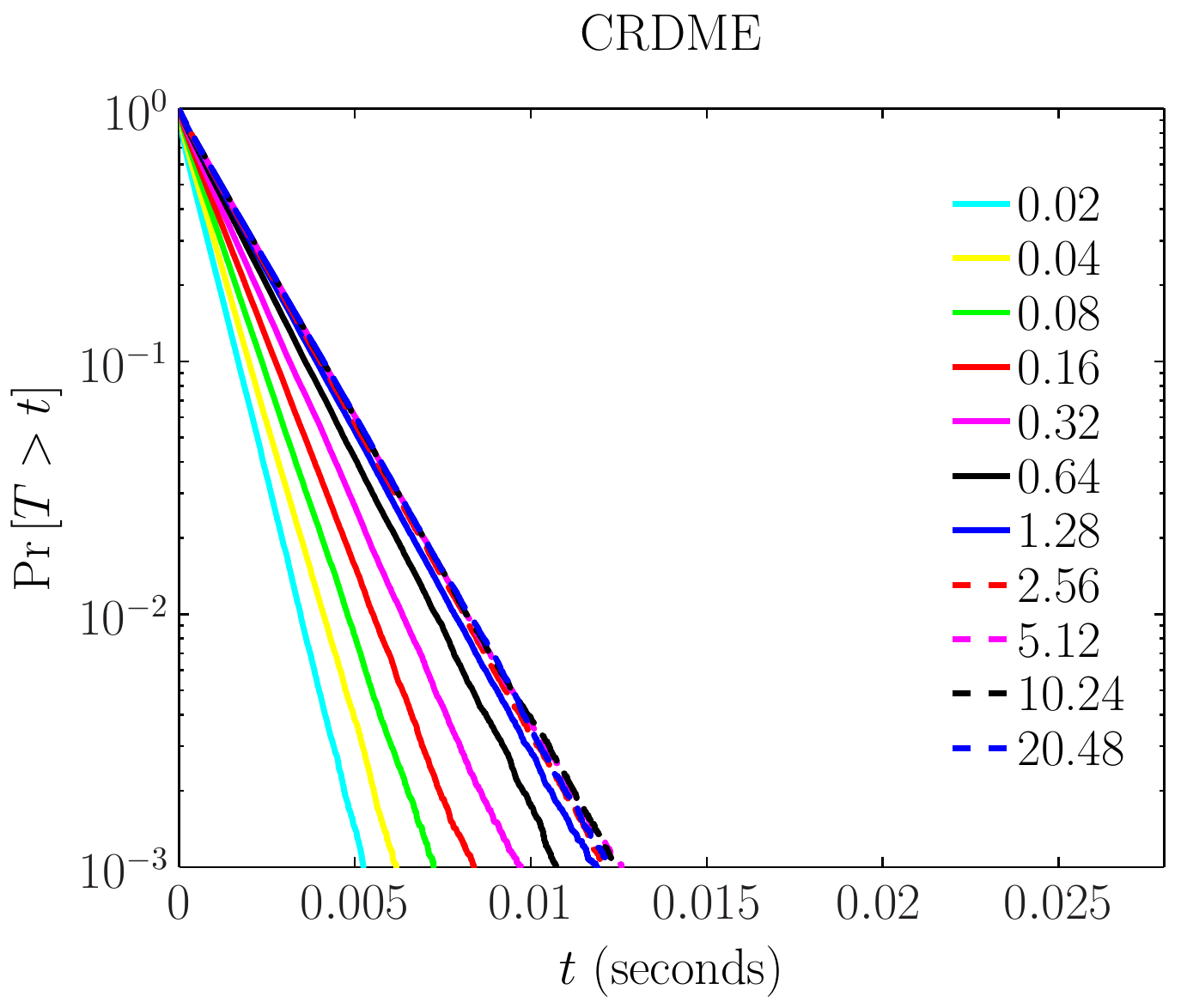}}  
  }
  \vspace{-20pt}
  \caption{ \small Survival time distributions vs. $t$ for the CRDME and RDME
    when $\lambda = 10^9 \, \textrm{s}^{-1}$. Each curve was estimated
    from $128000$ simulations. The legends give the ratio, $\rb / h$,
    used for each curve (note that $\rb = 1 \textrm{nm}$ was fixed and
    $h$ successively halved). We see the convergence of the survival
    time distributions for the CRDME (up to sampling error), while the
    survival time distributions in the RDME diverge. }
  \label{fig:survTimePlot}
\end{figure*}
\begin{figure*}
  \centering
  \scalebox{.54}{\includegraphics{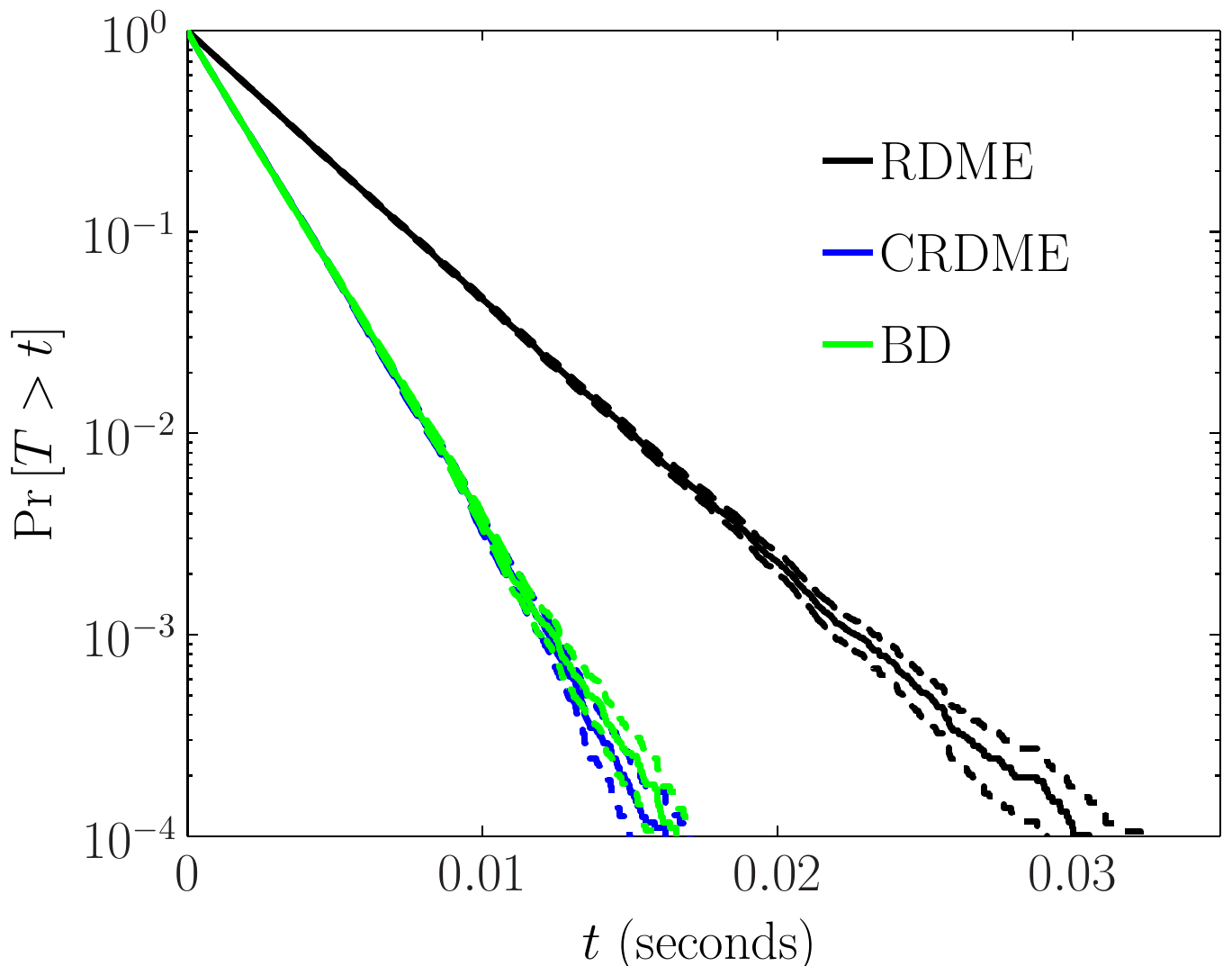}}  
  \caption{ \small \comment{Survival time distributions vs. $t$ for
      the RDME, CRDME, and from Brownian dynamics simulations (BD)
      when $\lambda = 10^9 \,\text{s}^{-1}$. The RDME and CRDME curves
      correspond to the survival time distributions shown in
      Fig.~\ref{fig:survTimePlot} for the largest value of $\rb/h$.
      For each curve $95\%$ confidence intervals are drawn with dashed
      lines (in the same color as the corresponding survival time
      distribution). They were determined using the Matlab
      \texttt{ecdf} routine. Each curve was estimated from $128000$
      simulations. To statistical error the CRDME and BD simulations
      agree, demonstrating that the CRDME has recovered the survival
      time distribution of the Doi model.}  }
  \label{fig:survTimePlotBD}
\end{figure*}
\begin{figure*}[t]
  \centering
  \subfloat[]{
    \scalebox{.54}{\includegraphics{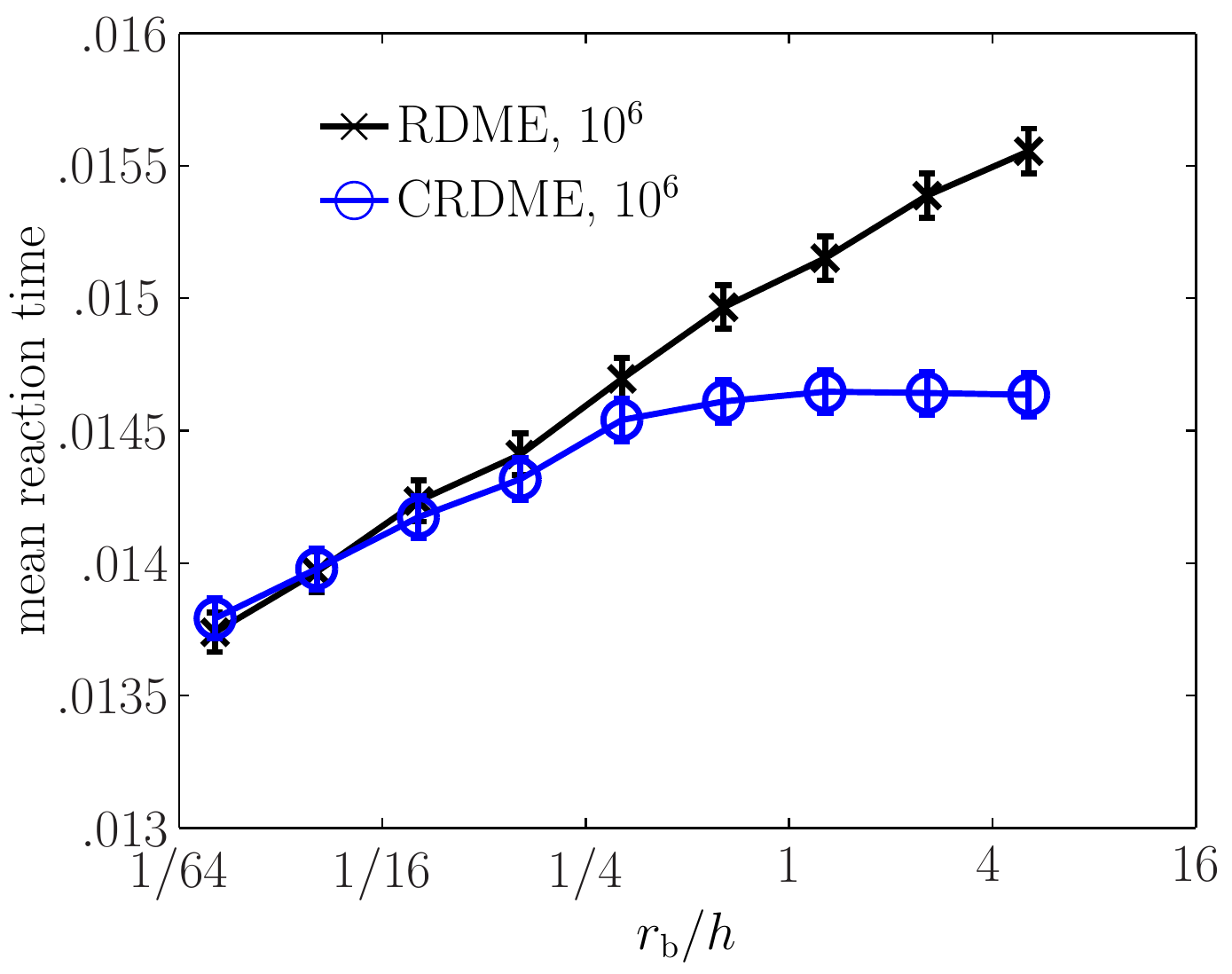}}  
  } \hspace{10pt}   
  \subfloat[]{
    \scalebox{.54}{\includegraphics{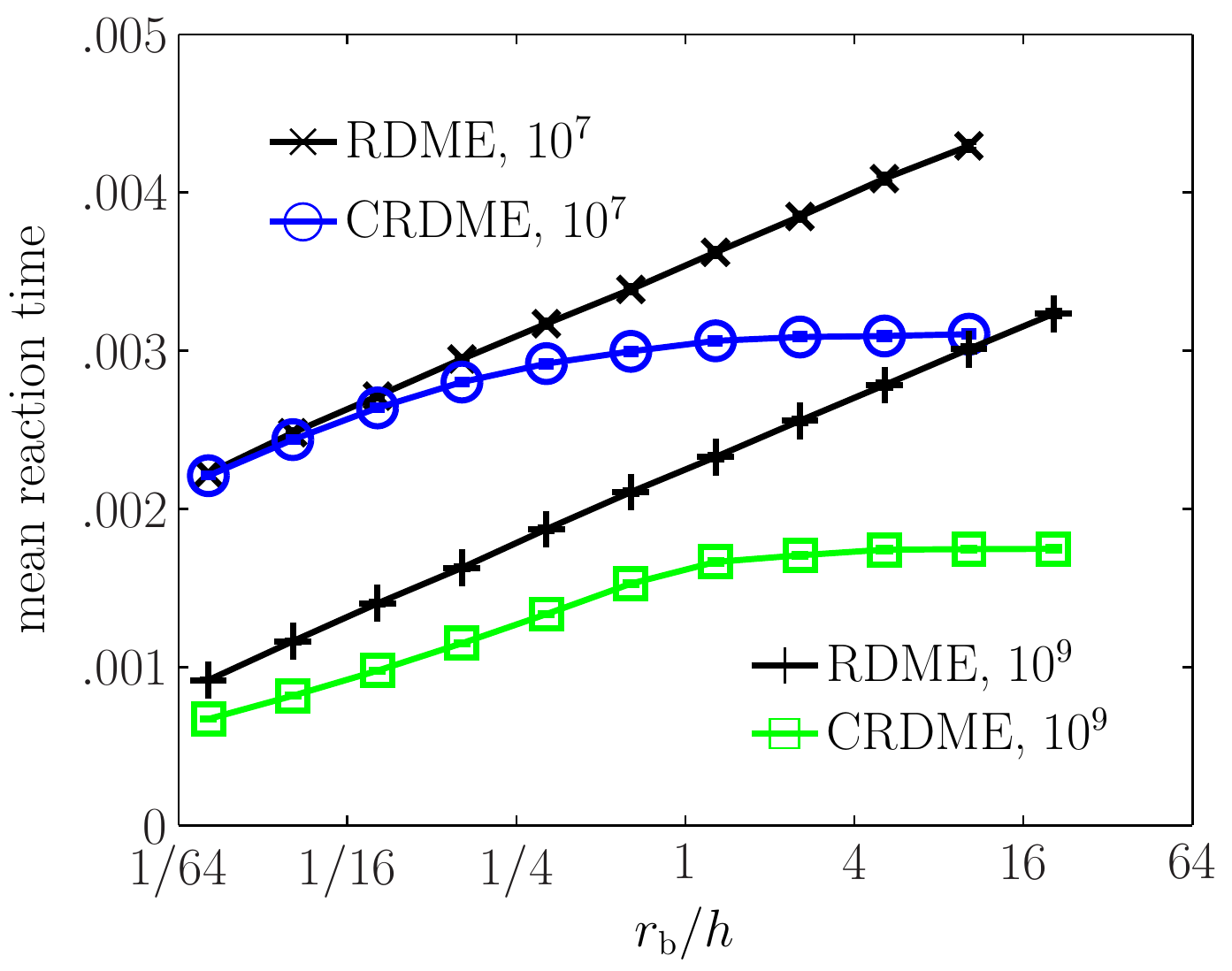}}  
  }
  \vspace{-20pt}
  \caption{  \small
    Mean reaction times vs. $\rb / h$ as $h$ decreases by
    factors of two. Legends give the value of $\lambda$ for each curve
    (with units of $\textrm{s}^{-1}$). Each mean reaction time was
    estimated from $128000$ simulations. Note that $95\%$ confidence
    intervals are drawn on each data point. (For some points they are
    smaller than the marker labeling the point.) Since we use a
    logarithmic $x$-axis, we see that for $h$ sufficiently small the
    mean reaction time in the standard RDME diverges like $\ln(h)$,
    while in the CRDME the mean reaction time converges to a finite
    value.  }
  \label{fig:meanBTPlot}
\end{figure*}
Let $T$ denote the random variable for the time at which the two
molecules react. The survival time distribution, $\prob \brac{T > t}$,
is
\begin{equation*}
  \prob \brac{T > t} = \int_0^L \int_0^L p(\vx,\vy,t) \, d \vx \, d \vy.
\end{equation*}
Note that the reaction time distribution, $\prob \brac{T < t} = 1 -
\prob \brac{T > t}$. We estimate $\prob{\brac{T > t}}$ from the
numerically sampled reaction times using the MATLAB \texttt{ecdf}
function.  In Fig.~\ref{fig:survTimePlot} we show the convergence (to
within sampling error) of the estimated survival time distributions
for the CRDME (right figure) as $h \to 0$ (for $\lambda = 10^9
\textrm{s}^{-1}$).  In the left figure we show the divergence as $h
\to 0$ of the estimated survival time distributions in the RDME. The
continuing rightward shift of the distribution as the mesh width is
decreased to twenty times finer than the reaction radius shows the
divergence of the reaction time to infinity. 

\comment{To confirm that the CRDME was converging to the solution of
  the Doi model we repeated these studies for $\lambda = 10^{9} \,
  \textrm{s}^{-1}$ using the Brownian dynamics (BD) method
  of~\cite{ErbanChapman2009,ErbanChapman2011}. In contrast to the RDME
  and CRDME, BD methods approximate the stochastic process describing
  the Brownian motion and reaction of the two molecules by
  discretization in time (instead of in space). For all BD simulations
  we used a fixed timestep, $dt = 10^{-10} \, \textrm{s}$. We refer
  the reader to~\cite{ErbanChapman2011} for details of the specific BD
  method we used. Fig.~\ref{fig:survTimePlotBD} illustrates that the
  survival time distribution in the CRDME and BD method agree to
  statistical error (using the CRDME with the smallest value of $h$
  from Fig~\ref{fig:survTimePlot}). The CRDME therefore recovers the
  reaction time statistics of the Doi model as $h \to 0$.}

\begin{figure}[t]
  \centering
  \scalebox{.54}{\includegraphics{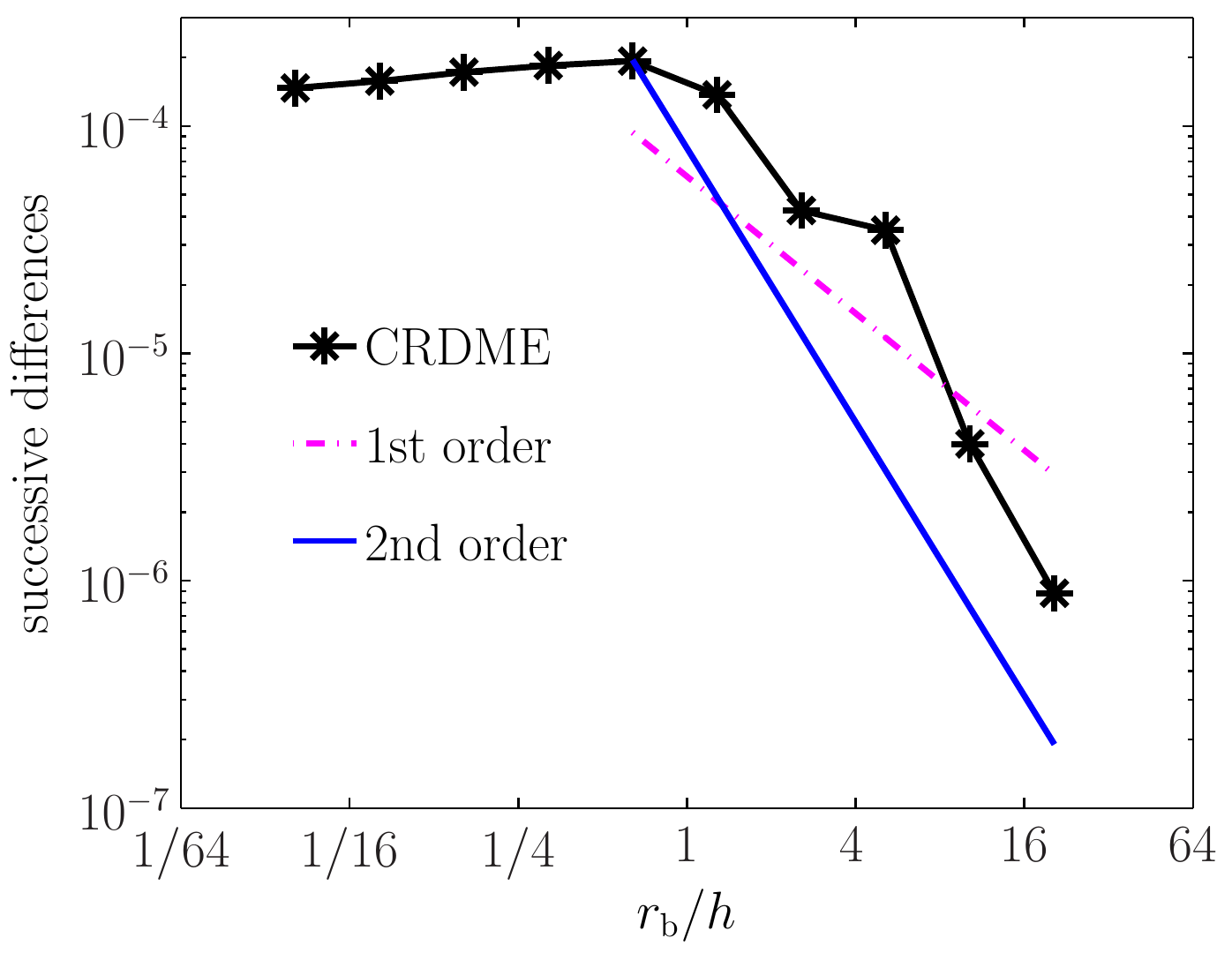}}  
  \caption{ \small Difference between successive points on the $\lambda =
    10^9$ CRDME curve in Fig.~\ref{fig:meanBTPlot} vs.  $\rb / h$. The
    smaller of the two $h$ values is used to label each point. The
    first and second order curves scale like $h$ and $h^2$
    respectively. Observe that the effective convergence rate to zero
    is closer to $O(h^2)$ than $O(h)$. }
  \label{fig:meanBTErr}
\end{figure}
The mean reaction time, $\avg{T}$, is given by
\begin{equation*}
  \avg{T} = \int_{0}^{\infty} \prob \brac{T > t} \, dt. 
\end{equation*}
We estimated $\avg{T}$ from the numerically sampled reaction times by
calculating the sample mean.  Fig.~\ref{fig:meanBTPlot} shows the
estimated mean reaction times for the CRDME and RDME models as
$\lambda$ and $\rb /h$ are varied (note the $x$-axis is logarithmic in
$\rb / h$). We see that as $h \to 0$ the sampled mean reaction times
in the CRDME converge to a fixed value. \comment{For $\lambda = 10^{9}
  \, \text{s}^{-1}$, the mean reaction time in the CRDME with the
  finest $h$ value, $0.0017471 \, \textrm{s}$, agreed with that found
  from the BD simulations used in Fig~\ref{fig:survTimePlotBD},
  $0.0017481 \, \textrm{s}$, to statistical error (\textit{i.e.}
  within $95\%$ confidence intervals, slightly less than $\pm 10^{-5}
  \, \textrm{s}$ about each mean value).}  The rate of convergence in
the CRDME for $\lambda = 10^9$ is illustrated in
Fig.~\ref{fig:meanBTErr}.  There we plot the difference between
successive estimated mean reaction times as $h$ is halved.  For small
values of $h$ this difference is seen to converge close to second
order (as illustrated by the slope of the solid blue line).

In contrast, Fig.~\ref{fig:meanBTPlot} shows that the sampled mean
reaction time in the RDME diverges like $\ln(h)$ as discussed
in~\cite{IsaacsonRDMELims,Hellander:2012jk}. For all $\lambda$ values
the sampled mean reaction time in the RDME converges to that of the
CRDME as $\rb / h \to 0$. As $\lambda$ is decreased we see agreement
between the RDME and CRDME for a larger range of $\rb / h$ values.

Figs.~\ref{fig:survTimePlot}, \ref{fig:meanBTPlot},
and~\ref{fig:meanBTErr} demonstrate that, in contrast to the RDME, the
reaction time statistics in the CRDME converge as $h \to 0$.
\comment{For $\lambda = 10^{9} \, \text{s}^{-1}$ we have verified the
  reaction time statistics of the CRDME recover those of the Doi model
  by comparison with BD simulations.} In the large lattice limit that
$\rb / h \to 0$ we see that reaction time statistics of the RDME
converge to those in the CRDME.  Hence we may interpret the RDME as an
approximation to the CRDME for $\rb / h \ll 1$.  The \emph{accuracy}
of using this approximation to describe the reaction-diffusion process
in the Doi model will then depend on the relative sizes of $\lambda$,
$D$, and $\rb$ as we discuss in the next section (and illustrated in
Fig.~\ref{fig:meanBTPlot}).

\section{RDME as an approximation of the CRDME for $\rb/h \ll
    1$} \label{S:rdmeApproxCRDME} 
\begin{figure}[t] 
  \centering
  \scalebox{.54}{\includegraphics{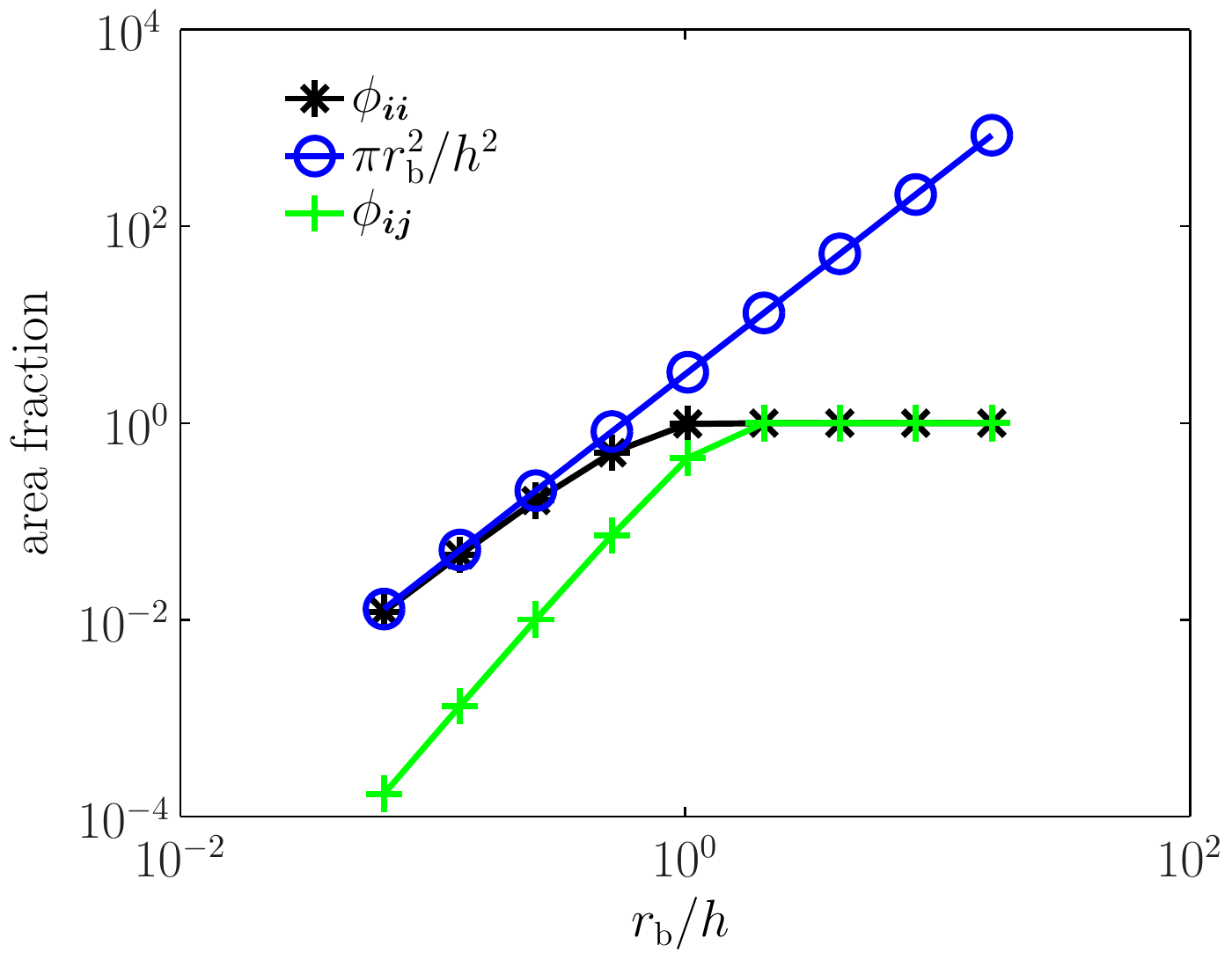}}
  \caption{ \small Comparison of area fraction when both molecules are in the
    same square, $\phi_{\bi \bi}$, with area fraction when the two
    molecules are in neighboring squares, $\phi_{\bi \bj}$. We see
    that the area fraction decreases faster as $\rb / h \to 0$ when
    the two molecules are in different squares than when they are in
    the same square.  Moreover, as $\rb / h \to 0$ the area fraction
    $\phi_{\bi \bi}$ approaches $\pi \rb^2 / h^2$ as derived
    in~\eqref{eq:RDMERxRate}.  }
  \label{fig:centVoxVol} 
\end{figure}
We now show that the RDME may be interpreted as an asymptotic
approximation to the CRDME for $\rb / h \ll 1$, with the accuracy of
this approximation depending on the size of $D$, $\rb$, and $\lambda$.
In the standard RDME two molecules can only react when within the same
$d$-dimensional voxel.  (If $k$ denotes a macroscopic bimolecular
reaction rate, the probability per unit time the molecules react is
usually chosen to be $k / h^d$, see~\eqref{eq:RDMEPartEqRxOp}.)  In
contrast, our new model allows two molecules to react when in
neighboring voxels.  Even for \emph{large} values of $h$, the volume
fraction $\phi_{\bi \bj}$ will be non-zero when $\bi$ and $\bj$ are
neighboring voxels.  That said, for $\bj \neq \bi$ the volume fraction
$\phi_{\bi \bj}$ will approach zero quicker as $\rb / h \to 0$ than
$\phi_{\bi \bi}$. This relationship is shown in two-dimensions ($d=2$)
in Fig.~\ref{fig:centVoxVol}. (We describe how the area fractions were
calculated in Appendix~\ref{S:reactTransitionRates}.)

We therefore expect that, asymptotically, when $\rb / h \to 0$ the
particles will effectively only react when $\bj = \bi$. In this case
\begin{align*}
  \abs{\mathcal{R} \cap V_{\bi \bi}} &= 
  \int_{\mathcal{R} \cap V_{\bi\bi}} \, d\vy \, d\vx \\
  &\approx \int_{\brac{-h/2,h/2}^d} \int_{\{\vy \mid \abs{\vx - \vy} < \rb \} } \, d\vy \, d\vx\\
  &= \abs{B_{\rb}} h^d,
\end{align*}
where $\abs{B_{\rb}}$ denotes the volume of the $d$-dimensional sphere
of radius $\rb$.  The reaction rate when both particles are at the
same position, $\bj = \bi$, is then
\begin{equation} \label{eq:RDMERxRate}
  \lambda \phi_{\bi \bi} = \frac{\lambda \abs{\mathcal{R} \cap V_{\bi \bi}}}{\abs{V_{\bi \bi}}} \approx
  \frac{\lambda \abs{B_{\rb}}}{h^d}.
\end{equation}
This corresponds to the choice of bimolecular reaction rate $k =
\lambda \abs{B_{\rb}}$ in the standard RDME. With this choice, when
$\rb / h \to 0$ the RDME may be interpreted as an asymptotic
approximation of the CRDME. This approximation is illustrated in
Fig.~\ref{fig:meanBTPlot}.

In three-dimensions, $d = 3$, the macroscopic bimolecular reaction
rate $k = \lambda \abs{B_{\rb}}$ also arises as the leading order
asymptotic expansion as $\rb \to 0$, $\lambda \to 0$, or $D = \DA +
\DB \to \infty$ of the diffusion limited bimolecular reaction rate for
the Doi model~\eqref{eq:DoiSimpEq}, $k_{\textrm{Doi}}$.
In~\cite{ErbanChapman2009} the latter was found to be
\begin{equation*}
  k_{\textrm{Doi}} = 4 \pi D \rb 
  \paren{ 1 - \frac{1}{\rb}\sqrt{\frac{D}{\lambda}} \tanh \paren{\rb \sqrt{\frac{\lambda}{D}}}}.
\end{equation*}
(Note, the more well-known Smoluchowski diffusion limited reaction
rate~\cite{SmoluchowskiDiffLimRx,KeizerJPhysChem82},
$k_{\textrm{Smol}} = 4 \pi D \rb$, is recovered in the limit $\lambda
\to \infty$.) As $\rb \sqrt{\lambda/D} \to 0$, $k_{\textrm{Doi}} \sim
\lambda (4 \pi \rb^3 / 3) = \lambda \abs{B_{\rb}}$.  We thus have that
the CRDME recovers this well-mixed reaction rate as $\rb / h \to 0$.
The smaller $\rb \sqrt{\lambda/D}$, the better the RDME should
approximate the CRDME for fixed $\rb / h$. 

When modeling three-dimensional biological systems, if $\rb
\sqrt{\lambda / D}$ is sufficiently small $h$ may simply be chosen to
accurately model molecular diffusion by a continuous-time random walk.
In this case, if $\rb / h \ll 1$ we may approximate the CRDME by the
standard RDME. When these assumptions break down we need to decrease
$h$ and incorporate reactions between molecules in neighboring voxels
with reactive transition rates $\lambda \phi_{\bi \bj}$.  The
reaction and diffusion processes could potentially be decoupled by
choosing separate meshes for each (as in~\cite{ElfPNASRates2010}).
That said, this process should be done so as to provide a convergent
approximation of~\eqref{eq:DoiSimpEq}.

\section{Reaction Product Placement} \label{S:crdmeRxProd}
Consider again the bimolecular reaction, $\textrm{A} + \textrm{B} \to
\textrm{C}$. As previously described, the fundamental problem that
prevents the standard RDME from converging to a reasonable continuous
particle method is the inability for reactants to find each other as
$h \to 0$. We therefore expect that any reasonable approach to placing
a newly created $\textrm{C}$ molecule will not impede the convergence
of the CRDME. In this section we propose one possible scheme for use
in the CRDME, consistent with the Doi model of
Section~\ref{S:mathModelFormulation}.  In the Doi model a newly
created $\textrm{C}$ molecule was placed at the center of the line
connecting the reacting $\textrm{A}$ and $\textrm{B}$ molecules. As
molecules are assumed uniformly distributed within voxels in the
CRDME, the method we now propose randomly places a newly created
molecule in one of several possible voxels surrounding those
containing two reactants.

To illustrate our approach we restrict to the two-particle $\textrm{A}
+ \textrm{B} \to \textrm{C}$ reaction for molecules in $\R^d$. As in
Section~\ref{S:CRDME} we denote by $p(\vx,\vy,t)$ the probability
density in the Doi model that the molecules of species $\textrm{A}$
and $\textrm{B}$ have not reacted and are located at $\vx$ and $\vy$
respectively at time $t$.  $p(\vx,\vy,t)$ still
satisfies~\eqref{eq:DoiSimpEq}, and will have the CRDME
approximation~\eqref{eq:RDMENewEq}. Let $u(\vq,t)$ denote the
probability density that a $\textrm{C}$ molecule has been created and
is located at position $\vq \in \R^{d}$ at time $t$.  If $\DC$ denotes
the diffusion constant of the $\textrm{C}$ molecule then, similar to
the in-flux term in the general Doi reaction
operator~\eqref{eq:doiRxOp}, we find $u(\vq,t)$ satisfies
\begin{equation} \label{eq:DoiSimpEqUnbind}
  \begin{aligned}
  \PD{u}{t} &= \DC \lap_{\vq} u + \lambda \int_{\R^d} \int_{\R^d} 
  \delta \paren{ \frac{\vx + \vy}{2} - \vq} \ind_{\mathcal{R}} \paren{\abs{\vx - \vy}}
  p(\vx,\vy,t) \, d \vx \, d \vy \\
    &= \DC \lap_{\vq} u + 2^d \lambda \int_{B_{\frac{\rb}{2}}(\vq)}
     p(\vx, 2 \vq - \vx, t) \, d \vx. 
  \end{aligned}
\end{equation}
Here $\lap_{\vq}$ denotes the Laplacian in $\vq$, while the second
term corresponds to the in-flux of probability density produced by the
reaction of the $\textrm{A}$ and $\textrm{B}$ molecules.
$B_{\frac{\rb}{2}}(\vq)$ denotes the hypersphere of radius $\rb/2$
centered at $\vq$.

We now derive a master equation approximation
of~\eqref{eq:DoiSimpEqUnbind}, consistent with the CRDME
approximation~\eqref{eq:RDMENewEq} of~\eqref{eq:DoiSimpEq}.
Discretize $\R^d$ into a lattice of cubic voxels of length $h$ indexed
by $\bk \in \Z^d$. The $\bk$th voxel is labeled by $V_{\bk}$, with
$\vq_{\bk} = \bk h$ denoting the center of the voxel. We make the
approximation that the probability the $\textrm{C}$ molecule exists
and is located in the voxel $V_{\bk}$ at time $t$ is given by
$U_{\bk}(t) = u(\vq_{\bk},t) \abs{V_{\bk}}$. Using this assumption we
construct a finite volume discretization of~\eqref{eq:DoiSimpEqUnbind}
by integrating both sides of~\eqref{eq:DoiSimpEqUnbind} over
$V_{\bk}$.  We use the same finite volume approximation of the
Laplacian as before, obtaining a discretized Laplacian in the $\vq$
coordinate as defined in~\eqref{eq:lapDef}. Denote this discrete
Laplacian by $\DC L_h^{\textrm{C}}$.  The integral of the incoming
flux term in~\eqref{eq:DoiSimpEqUnbind} is approximated by
\begin{align*}
\int_{V_{\bk}} \int_{\R^d} \int_{\R^d} 
\delta \paren{ \frac{\vx + \vy}{2} - \vq} \ind_{\mathcal{R}} &\paren{\abs{\vx - \vy}}
p(\vx,\vy,t) \, d \vx \, d \vy \, d \vq \\
  &= \int_{\R^{d}} \int_{\R^{d}} \ind_{V_k} \paren{\frac{\vx + \vy}{2}} 
  \ind_{\mathcal{R}} \paren{\abs{\vx-\vy}} p(\vx,\vy,t) \, d \vx \, d \vy \\
  &\approx \sum_{\substack{\bi \in \Z^d\\ \bj \in \Z^d}} \gamma_{\bi \bj}^{\bk} \phi_{\bi \bj} P_{\bi,\bj}(t).
\end{align*}

Here $\ind_{V_{\bk}}(\vq)$ denotes the indicator function of the set
$V_{\bk}$ and $P_{\bi,\bj}(t)$ denotes the solution to the
CRDME~\eqref{eq:RDMENewEq}.  $\gamma_{\bi \bj}^{\bk}$ is defined by
\begin{equation} \label{eq:crdmePlacementProb}
  \gamma_{\bi \bj}^{\bk} = \begin{cases}
    \frac{1}{\phi_{\bi \bj} \abs{V_{\bi \bj}}} \int_{V_{\bi \bj}} \ind_{V_k} \paren{\frac{\vx + \vy}{2}} 
    \ind_{\mathcal{R}} \paren{\abs{\vx-\vy}} d \vx \, d \vy, &\phi_{\bi \bj} \neq 0,\\
    0, & \phi_{\bi \bj} = 0,
    \end{cases}
\end{equation}
and gives the probability that when $\vx \in V_{\bi}$ and $\vy \in
V_{\bj}$ react the resultant $\textrm{C}$ molecule is placed in voxel
$V_{\bk}$. Note,
\begin{equation*}
  \sum_{\bk \in \Z^d} \gamma_{\bi \bj}^{\bk} = 1
\end{equation*}
when $\phi_{\bi \bj} \neq 0$.  In practice, it is possible to
calculate $\gamma_{\bi \bj}^{\bk}$ by modification of the algorithm
given in Appendix~\ref{S:reactTransitionRates} for calculating
$\phi_{\bi \bj}$.

We therefore arrive at the following master equation approximation
of~\eqref{eq:DoiSimpEqUnbind}
\begin{equation} \label{eq:CRDMESimpEqUnbind}
  \D{U_{\bk}}{t} = \DC L_h^{\textrm{C}} U_{\bk} + \lambda \sum_{\substack{\bi \in \Z^d\\ \bj \in \Z^d}} \gamma_{\bi \bj}^{\bk} \phi_{\bi \bj} P_{\bi,\bj}(t).
\end{equation}
Using the CRDME model given by~\eqref{eq:RDMENewEq}
and~\eqref{eq:CRDMESimpEqUnbind} we are then led to a CRDME
approximation of the general Doi model~\eqref{eq:doiFormEvolEq}.  Let
$\tilde{F}_{h}^{(a,b,c)}
\paren{\bja, \bjb, \bjc,t}$ denote the solution to the general CRDME,
analogous to the solution of the standard RDME, $F_{h}^{(a,b,c)}
\paren{\bja, \bjb, \bjc,t}$. $\tilde{F}_h$ also
satisfies~\eqref{eq:RDMEPartEq}, with the same diffusion operator,
$L_h$, but with the modified reaction operator, $\tilde{R}_h$, given
by
\begin{multline} \label{eq:genCRDMERxOp}
  \paren{ \tilde{R}_{h} \tilde{F}_{h}^{\paren{a,b,c}} } \beqparen{\bja,\bjb,\bjc,t} = 
  \lambda \bigg[  \sum_{l=1}^{c} \sum_{
    \substack{
      \bi \in \Z^d\\
      \bj \in \Z^d
    }}
  \gamma_{\bi \bj}^{\bjc_l}
  \phi_{\bi \bj} \tilde{F}_{h}^{\paren{a+1,b+1,c-1}}
  \beqparen{\bja \cup \bi, \bjb \cup \bj, \bjc \setminus \bjc_{l},t} \\
  - \sum_{l=1}^{a} \sum_{m=1}^{b} \phi_{\bja_{l} \bjb_{m}} \,
  \tilde{F}_{h}^{\paren{a,b,c}} \beqparen{\bja,\bjb,\bjc,t} \bigg].
\end{multline}

\section{Modified SSA  for Bimolecular Reactions Based on the CRDME} \label{S:crdmeSSA}
We conclude by summarizing how to modify the SSA to generate
realizations of the CRDME. Let $a_{\bi}$ denote the current number of
molecules of species $\textrm{A}$ in voxel $\bi$ in a simulation, with
$b_{\bj}$ and $c_{\bk}$ defined similarly.  We are then lead to the
following proposed modification of the SSA to handle the bimolecular
reaction $\textrm{A} + \textrm{B} \to \textrm{C}$ based on the CRDME
\begin{enumerate}
\item The probability per unit time a molecule of species $\textrm{A}$
  in voxel $\bi$ reacts with a molecule of species $\textrm{B}$ in
  voxel $\bj$ is given by the propensity $\lambda \phi_{\bi \bj}
  a_{\bi} b_{\bj}$.
\item Should such a reaction occur, update the system state so that
  \begin{enumerate}
  \item $a_{\bi} \from a_{\bi} - 1$.
  \item $b_{\bj} \from b_{\bj} - 1$.
  \item Randomly chose a voxel $\bk$ with probability $\gamma_{\bi \bj}^{\bk}$
    and update $c_{\bk} \from c_{\bk} + 1$. 
  \end{enumerate}
\item Recalculate any reaction or event times that depend on $a_{\bi}$,
  $b_{\bj}$, or $c_{\bk}$.
\end{enumerate}
For all diffusive transitions, zeroth order reactions, and first
order reactions the SSA remains the same as for the standard
RDME~\cite{ElfIEESys04}.

\section{Conclusion}
By discretizing the stochastic reaction-diffusion model of
Doi~\cite{DoiSecondQuantA,DoiSecondQuantB} we have derived a new
convergent reaction-diffusion master equation for $\textrm{A} +
\textrm{B} \to \textrm{C}$. We illustrated our discretization
procedure, and the convergence of the survival time distribution and
mean reaction time in the CRDME, for two molecules that undergo the
annihilation reaction $\textrm{A} + \textrm{B} \to \varnothing$. While
this special case is simplified compared to realistic biological
networks, it should be noted that the same reaction rates, $\lambda
\phi_{\bi \bj}$, are obtained by this discretization procedure for the
more general multiparticle Doi model.  This resulted in the general
CRDME for $\textrm{A} + \textrm{B} \to \textrm{C}$ given
by~\eqref{eq:doiFormEvolEq} with the reaction
operator~\eqref{eq:genCRDMERxOp}.  \comment{While we derived a
  convergent RDME by discretization of the Doi model in this work, we
  expect that a similar finite-volume discretization approach might
  also allow the derivation of a convergent RDME-like approximation to
  Smoluchowski models.}  




\begin{acknowledgments}
  SAI is supported by NSF grant DMS-0920886. SAI thanks I. Agbanusi, 
  D. Isaacson, and A. Steele for helpful comments and suggestions.
\end{acknowledgments}

\appendix

\section{Calculation of reaction transition rates} \label{S:reactTransitionRates}

Denote by $V_{\bi}$ the $d$-dimensional coordinate axis aligned
hypercube with sides of length $h$ centered at $\bi h$. With this
definition we may then write $V_{\bi \bj} = V_{\bi} \times V_{\bj}$.
We use $\hat{V}_{\bi}$ to denote this hypercube in the special case
that $h = 1$.  Finally, let $B_{\rb}(\vx)$ be the $d$-dimensional
hypersphere of radius $\rb$ about $\vx$. A convenient representation
for $\phi_{\bi \bj}$ we subsequently use is
\begin{align}
  \phi_{\bi \bj} &= \frac{1}{\abs{V_{\bi \bj}}} \int_{V_{\bi}}
  \int_{V_{\bj}} \ind_{R}(\abs{\vx - \vy}) \, d \vy \, d \vx \label{eq:phiDoubleIntDef} \\
  &= \frac{1}{\abs{V_{\bi \bj}}} \int_{V_{\bi}} \abs{B_{\rb}(\vx) \cap
    V_{\bj}} \, d \vx  \notag \\
  &= \int_{\hat{V}_{\vO}} \abs{ B_{\frac{\rb}{h}}(\vx) \cap \hat{V}_{\bj-\bi}} \, d \vx.
  \label{eq:phiIntDef}
\end{align}
(Here $\vO$ denotes the origin voxel.)  Hence we may interpret
$\phi_{\bi \bj}$ as the integral over the center of a hypersphere of
the volume of intersection between the hypersphere and a hypercube.
The final equation~\eqref{eq:phiIntDef} shows that $\phi_{\bi \bj}$
depends on only two quantities; the separation vector $\bj - \bi$ and
$\rb / h$. Also note that $\phi_{\bi \bj}$ will be zero once the
separation between all points in voxels $\bi$ and $\bj$ is more than
$\rb$. As such, in practice it is only necessary to calculate
$\phi_{\vO \bj}$ for a small number of voxels about the origin.

It is desirable to calculate $\phi_{\bi \bj}$ to near machine
precision to avoid the introduction of error from the use of incorrect
reactive transition rates. While this may seem an easy task, simply
calculating the hypervolume of intersection of $\mathcal{R}$ and
$V_{\bi \bj}$, it should be noted that these are four-dimensional
(six-dimensional) sets when the molecules are in two-dimensions
(three-dimensions).  Evaluating $\phi_{\bi \bj}$ by directly applying
quadrature to~\eqref{eq:phiDoubleIntDef} is complicated by the
discontinuous integrand. We have found that several standard
cubature~\cite{CubaHahn2005,CubatureMIT} and Monte Carlo
methods~\cite{CubaHahn2005} have difficultly evaluating such integrals
in reasonable amounts of computing time to high numerical precision
(absolute errors below $10^{-11}$). Since the
integral~\eqref{eq:phiIntDef} has a continuous integrand, which only
requires the intersection of two-dimensional (three-dimensional) sets
when the particles are each in two-dimensions (three-dimensions), we
focus on evaluating $\phi_{\bi \bj}$ through this representation.

To evaluate~\eqref{eq:phiIntDef} both efficiently and accurately it is
necessary to calculate the hypervolume of intersection given by the
integrand, $v_{\bj}(\vx) = \abs{B_{\frac{\rb}{h}}(\vx) \cap \hat{V}_{\bj}}$.
Our approach is based on writing this hypervolume as an integral and
then converting to a boundary integral through the use of the
divergence theorem.  That is,
\begin{align}
  v_{\bj}(\vx) &= \frac{1}{d} \int_{B_{\frac{\rb}{h}}(\vx) \cap \hat{V}_{\bj} } \nabla \cdot \vy \, d \vy, \notag \\
  &= \frac{1}{d} \int_{\partial(B_{\frac{\rb}{h}}(\vx) \cap \hat{V}_{\bj}) } \vy \cdot \veta(\vy) \, dS(\vy), \notag \\
  &= \frac{1}{d} \int_{\partial B_{\frac{\rb}{h}}(\vx)}  \paren{\vy \cdot \veta(\vy)} \ind_{\hat{V}_{\bj}}(\vy)\, dS(\vy) 
  +\frac{1}{d} \int_{\partial \hat{V}_{\bj}}  \paren{\vy \cdot \veta(\vy)} \ind_{B_{\frac{\rb}{h}}(\vx)}(\vy)\, dS(\vy). \label{eq:areaOfIsect}
\end{align}
Here $\partial M$ is used to denote the boundary of a manifold $M$,
$\veta(\vy)$ the outward normal to the boundary hypersurface at $\vy$,
and $dS(\vy)$ the hypersurface measure at $\vy$.

For simplicity, in the remainder we assume $d = 2$.  In this case we
have developed a fast method, requiring only a few minutes on a modern
laptop, that is able to evaluate~\eqref{eq:phiIntDef} to near machine
precision.  $v_{\bj}(\vx)$ is evaluated by calculating the
intersection points of the circle $\partial B_{\rb}(\vx)$ with the
square $\hat{V}_{\bj}$ numerically.  Once these points are known the
line integrals in~\eqref{eq:areaOfIsect} can be reduced to sums of
integrals over sub-arcs where the indicator function is identically
one or zero.  These integrals can be evaluated analytically.  Standard
adaptive numerical quadrature methods, such as the \texttt{dblquad}
routine in MATLAB, are then able to effectively integrate the area of
intersection function $v_{\bj}(\vx)$.  This method was used to
generate the area fractions in Fig.~\ref{fig:centVoxVol} and all SSA
simulations.


\bibliographystyle{aipnum4-1.bst}

%

\end{document}